\documentclass[11pt]{article}
\usepackage[british]{babel}
\usepackage{graphicx,float,multirow,fullpage,amsmath,amssymb}  

\usepackage{caption2}

\usepackage{color,hyperref, cite}
\usepackage{algorithm, algorithmic}
\usepackage{ulem}

\definecolor{darkblue}{rgb}{0.0,0.0,0.5}
\definecolor{darkgrey}{rgb}{0.5,0.1,0.1}
\hypersetup{colorlinks,breaklinks,
            linkcolor=darkblue,urlcolor=darkblue,
            anchorcolor=darkblue,citecolor=darkblue}

\newcommand{\eref}[1]{(\ref{#1})}

\DeclareTextSymbol{\degre}{T1}{6}
\DeclareTextSymbol{\degre}{OT1}{23}

\begin{document}
\begin{center}
{\bf{\LARGE DynPeak : An algorithm for pulse detection \\[0.2cm]
and frequency analysis in hormonal time series}}\\[0.5cm]
\end{center}

\begin{center}
{\bf Alexandre VIDAL}$^{1,\ast}$, 
{\bf Qinghua ZHANG}$^{2}$, 
{\bf Claire M\'EDIGUE}$^{3}$, \\
\vspace{0.1cm}
{\bf St\'ephane FABRE}$^{4}$,
{\bf Fr\'ed\'erique CL\'EMENT}$^{3}$
\end{center}
\vspace{-0.2cm}

\begin{flushleft}
\bf{1} Laboratoire Analyse et Probabilit\'es EA 2172, F\'ed\'eration de Math\'ematiques FR 3409, Universit\'e d'\'Evry-Val-d'Essonne, 23 Boulevard de France, 91037, Evry, France.
\\[0.2cm]
\bf{2} Project-Team SISYPHE, INRIA Rennes-Bretagne Atlantique Research Centre, \\
Campus de Beaulieu, 35042, Rennes Cedex, France.
\\[0.2cm]
\bf{3} Project-Team SISYPHE, INRIA Paris-Rocquencourt Research Centre, \\
Domaine de Voluceau Rocquencourt - B.P. 105, 78153 Le Chesnay, France.
\\[0.2cm]
\bf{4} Laboratoire de Physiologie de la Reproduction et des Comportements, UMR 85 INRA, UMR 6175 CNRS, INRA Tours Research Center - Universit\'e F. Rabelais de Tours - IFCE, 37380 Nouzilly, France.
\\[0.2cm]
$\ast$ E-mail: \href{mailto:alexandre.vidal@univ-evry.fr}{\tt alexandre.vidal@univ-evry.fr}
\end{flushleft}

\section*{Abstract}

The endocrine control of the reproductive function is often studied from the analysis of luteinizing hormone (LH) pulsatile secretion by the pituitary gland. Whereas measurements in the cavernous sinus cumulate anatomical and technical difficulties, LH levels can be easily assessed from jugular blood. However, plasma levels result from a convolution process due to clearance effects when LH enters the general circulation. Simultaneous measurements comparing LH levels in the cavernous sinus and jugular blood have revealed clear differences in the pulse shape, the amplitude and the baseline. Besides, experimental sampling occurs at a relatively low frequency (typically every 10 min) with respect to LH highest frequency release (one pulse per hour) and the resulting LH measurements are noised by both experimental and assay errors. As a result, the pattern of plasma LH may be not so clearly pulsatile. Yet, reliable information on the InterPulse Intervals (IPI) is a prerequisite to study precisely the steroid feedback exerted on the pituitary level. Hence, there is a real need for robust IPI detection algorithms.

In this article, we present an algorithm for the monitoring of LH pulse frequency, basing ourselves both on the available endocrinological knowledge on LH pulse (shape and duration with respect to the frequency regime) and synthetic LH data generated by a simple model. We make use of synthetic data to make clear some basic notions underlying our algorithmic choices. We focus on explaining how the process of sampling affects drastically the original pattern of secretion, and especially the amplitude of the detectable pulses. We then describe the algorithm in details and perform it on different sets of both synthetic and experimental LH time series. We further comment on how to diagnose possible outliers from the series of IPIs which is the main output of the algorithm.

\section*{Introduction}

The neuroendocrine axes play a major part in controlling the main physiological functions (me\-tab\-o\-lism, growth, development and reproduction). The connection between the central nervous system and the endocrine system takes place on the level of the hypothalamus, where endocrine neurons are able to secrete hormones that target the pituitary gland.  In birds and mammals, a dedicated portal system (the pituitary portal system) joins the hypothalamus and pituitary gland together. The anterior lobe of the pituitary gland (adenohypophysis) produces different hormones, which target either other endocrine glands (releasing their hormones directly into the bloodstream), exocrine glands (releasing their hormones into dedicated ducts) or non-secreting organs.

We will be particularly interested in the gonadotropic axis, that is named according to its most downstream component, the gonads (ovaries in females, testes in males). The reproductive axis is under the control of the gonadotropin-releasing hormone (GnRH), which is secreted in pulses from specific hypothalamic areas. GnRH effects on its target cells depend critically on pulse frequency and ultimately result in the differential secretion patterns of the luteinizing hormone (LH) and follicle-stimulating hormone (FSH). LH secretion pattern is clearly pulsatile, while FSH pattern is not. LH and FSH control the development of germinal cells within the gonads and the secretory activity of somatic cells. In turn, hormones secreted by the gonads (steroid hormones such as androgens, progestagens and oestrogens or peptidic hormones such as inhibin) modulate the secretion of GnRH, LH and FSH within intertwined feedback loops.

Whereas measurements of GnRH levels (in either the pituitary portal blood or the cerebrospinal fluid) cumulate anatomical and technical difficulties, LH levels can be easily assessed from jugular blood. In females, there is a clear modulation of LH pulse frequency along an ovarian cycle \cite{Sollenberger_90}. Pulse frequency is much lower in the luteal, progesterone-dominated phase compared to the follicular, oestradiol-dominated phase.  Apart from the period surrounding ovulation, there is a good correlation between GnRH and LH pulses \cite{Moenter_90, Moenter_92}, so that a precise determination of LH pulse frequency is valuable to investigate the feedback effects of gonadal hormones in different physiological or pathological situations.

LH plasma levels result from a convolution process. The instantaneous LH release rate from the pituitary gland is pulsatile, but as soon as LH enters the general circulation, it is subject to clearance  effects.  Simultaneous measurements  of LH levels in the cavernous sinus and jugular blood \cite{Clarke_02} have revealed clear differences in the pulse shape and amplitude as well as in the baseline. Besides, experimental sampling occurs at a relatively low frequency (typically every $10$ min, \cite{Drouilhet_10, Baird_81a, Baird_81b}) with respect to LH highest frequency release (one pulse per hour) and the resulting LH measurements are noised by both experimental and assay errors. As a result, the pattern of plasma LH may be not so clearly pulsatile. Yet, reliable information on the interpulse intervals (IPI) is a prerequisite to study precisely the steroid feedback exerted on the pituitary level.  Hence, there is a real need for robust IPI detection algorithms.

In this article, we present an algorithm for the monitoring of LH pulse frequency, basing ourselves both on the available endocrinological knowledge on LH pulse (shape and duration with respect to the frequency regime) and synthetic LH data generated by a simple model. We make use of synthetic data to make clear some basic notions underlying our algorithmic choices. We focus on explaining how the process of sampling affects drastically the original pattern of secretion, and especially the amplitude of the detectable pulses. We then describe the algorithm in details and perform it on different sets of both synthetic and experimental LH time series. We further comment on how to diagnose possible outliers from the series of IPIs which is the main output of the algorithm.

\section*{Methods}

\subsection*{A mathematical generator of synthetic LH time series}

\begin{figure}[H]
\centering
\includegraphics[scale=0.4]{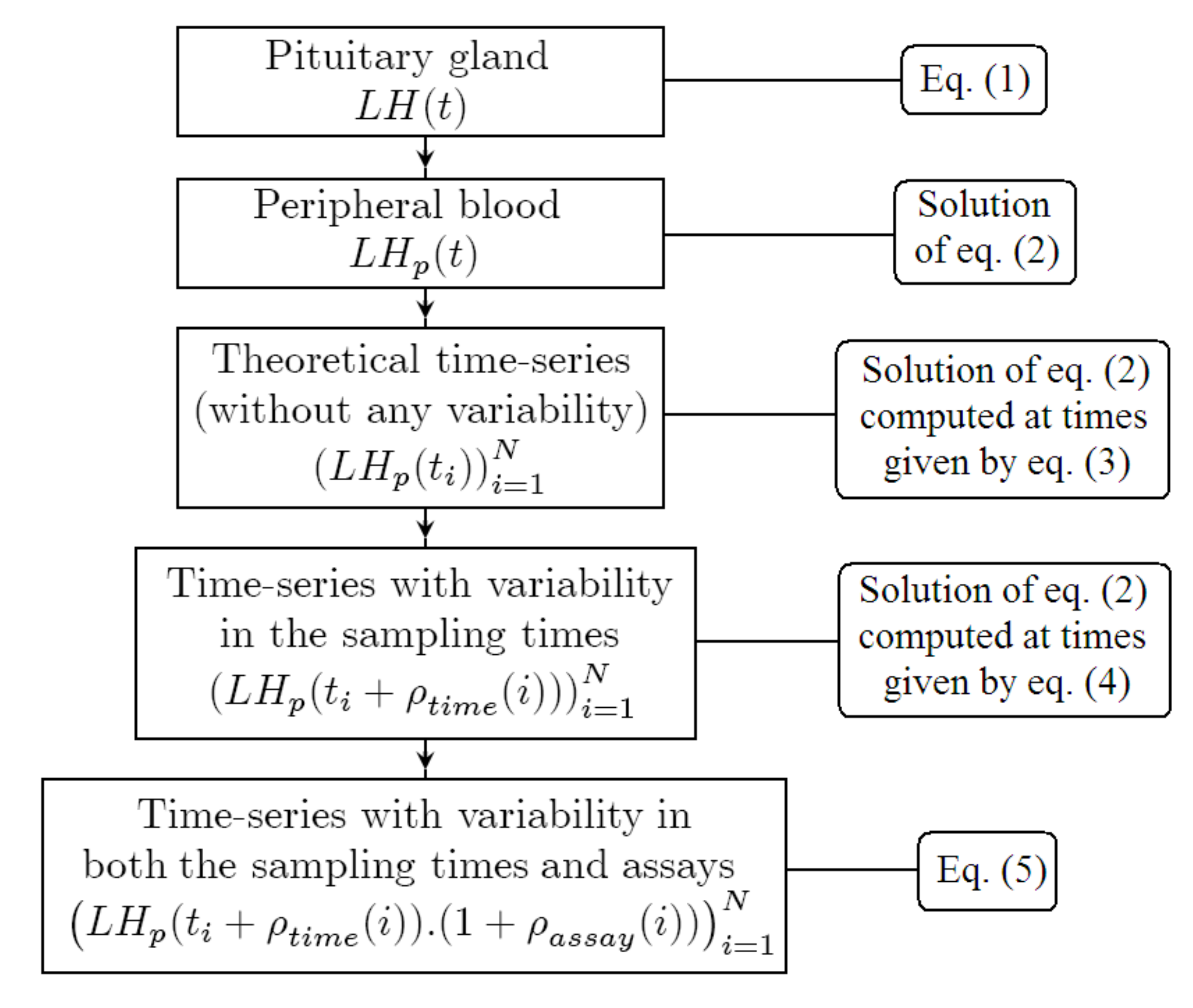}
\caption{\bf Diagram of the synthetic sampling process.}
\label{DiagrammeSampling}
\end{figure}

Basing ourselves on a simple model of plasma LH level introduced in \cite{Vidal_09}, we illustrate the effects of the sampling process upon a LH signal. This model combines a representative function of the pulsatile LH secretion by the pituitary gland with a term accounting for the clearance from the blood. The synthetic sampling process is designed to reproduce as close to experiments as possible the variability in the sampling times and measurements. The different steps of this construction, as well as the links between the mathematical objects and what they represent in the biological context, are illustrated in the diagram of Figure \ref{DiagrammeSampling}.

\subsubsection*{Model of Luteinizing Hormone secretion}

The pituitary gland releases LH into the blood as successive spikes characterized by a quasi-instantaneous increase followed by slower (yet quite fast) decrease (see \cite{Clarke_02}). Hence, in our model, the LH release along a spike is approximated by a discontinuous function of time: the jump accounting for the instantaneous increase in the LH release is followed by a fast exponential decrease. The interspike interval is controlled by a function of time $P_{spike}(t)$ accounting for the varying release frequency. The spike amplitude is also subject to an inter-spike variability as well as to long-term changes partly due to the time variations in the stock of LH available for release. In our model, the amplitude is controlled by another function of time $M_{spike}(t)$. Hence, the instantaneous release of LH (expressed in ng/ml/min) in the blood by the pituitary gland is given by:
\begin{equation}
LH(t)=M_{spike}(t)\exp \left[ -k_{hl}\left( t-\left\lfloor \frac{t}{P_{spike}(t)}\right\rfloor P_{spike}(t)\right) \right]
\label{LHspike}
\end{equation}
where $\lfloor x \rfloor$ refers to the greatest integer smaller than $x$ (integer part). The $k_{hl}$ stiffness coefficient is directly linked to the spike half-life, $\tau _{hl}$, by the following relation:
\[
k_{hl}=\frac{\ln 2}{\tau _{hl}}
\]
We based our choice of parameter values on the few experiments that have investigated the LH release by the pituitary gland in the ewe from synchronous sampling in the jugular blood and cavernous sinus \cite{Clarke_02}. Accordingly, we chose a spike half-life $\tau_{hl}=20$ min (i.e. $k_{hl} \simeq 0.03466$ min$^{-1}$).

We represent the continuously measured LH blood level (expressed in ng/ml) as the solution of:
\begin{equation}
\frac{{\rm d}LH_{p}}{{\rm d}t}(t)=LH(t)-\alpha LH_{p}(t)
\label{LHplasma}
\end{equation}
where the LH release rate LH(t) is given by equation (\ref{LHspike}). Parameter $\alpha$ represents the instantaneous LH clearance rate from the blood. To be consistent with the one hour half-life of LH pulses in the jugular blood, we have fixed $\alpha =6$ min$^{-1}$.

\subsubsection*{Sampling protocol and assay variance}

To mimic the experimental protocol for LH data acquisition, we extract time series from the fine step numerical integration of equation (\ref{LHplasma}). This process is intended to obtain a time series of $N$ consecutive samples similar to experimental results, i.e. a finite sequence of $N$ couples $(t_i, A_i)$, with $i=1,2,...N$, where $A_i$ is the measured LH level at time $t_i$.

In most experiments, the LH data are retrieved at a fixed frequency. We describe below the corresponding synthetic process: the samples are obtained each $T_s$ minutes from the starting time $t=r$ min. Hence, the sampling times are 
\begin{equation} \label{Times_t_i}
t_i = r+T_s (i-1), \quad i=1,2,...N.
\end{equation}
Note that the first sample is retrieved at time $t_1=r$. Parameter $r \in [0,T_s]$ allows one to shift the beginning of the sampling process while using the same set of data simulated from equation (\ref{LHplasma}).

To take into account the inherent variability of the experimental sampling times, we compute times $\tau_i$ near the theoretical sampling time $t_i$:
\begin{equation} \label{Times_tau_i}
\tau_{i}=t_i+\rho_{time}(i), \quad i=1,2,...N.
\end{equation}
The random numbers $\rho_{time}(i)\in [-f,f]$, $f>0$, are generated from an almost uniform distribution (using the Mersenne Twister algorithm \cite{Matsumoto_98}). Then, we can compute the value $LH_p (\tau_i)$ of the solution of equation \eref{LHplasma} at each of these times.

\begin{figure}[!ht]
\centering
\includegraphics[width=15cm]{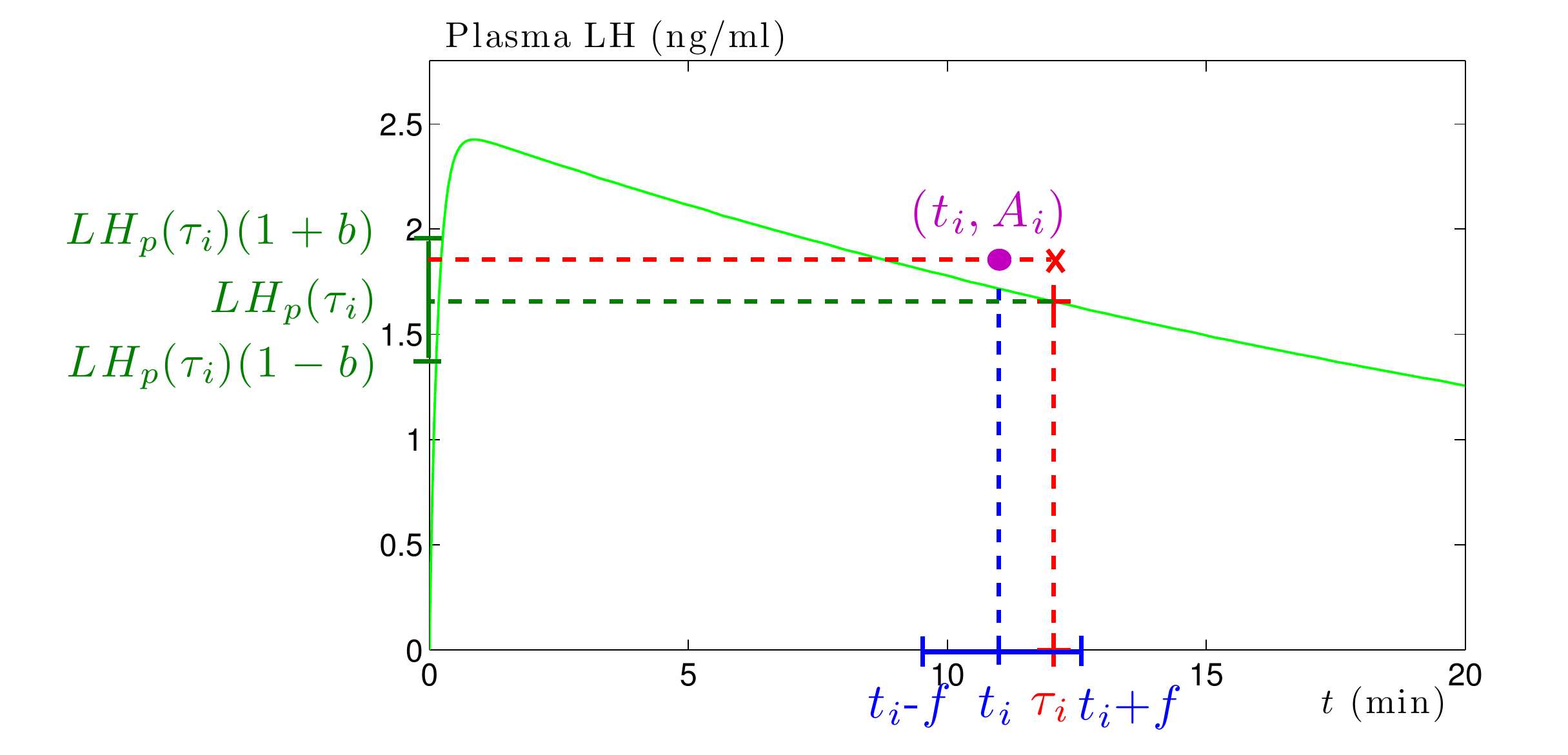}
\caption{
 {\bf  Computation of a synthetic sample.}
Solution $LH_p(t)$ of equation (\ref{LHplasma}) (green curve), retained sampling time $t_i$ (blue dotted line), real sampling time $\tau _i$ (red dotted line) randomly chosen in $[t_i-f, t_i+f]$ (blue interval), exact value $LH_p(\tau _i)$ of LH level at time $\tau_i$ (green dotted line), retrieved LH level $A_i$ (red dotted line) randomly chosen in $[LH_p(\tau_i)(1-b), LH_p(\tau_i)(1+b)]$ (green interval). The output of the $i$-th step of the sampling process is the couple $(t_i, A_i)$ (magenta disc) of time and corresponding LH level.}
\label{Construction_A_i}
\end{figure}

We also reproduce the LH assay variance by applying a multiplicative noise on each sample. We compute:
\begin{equation} \label{Samples}
A_i=LH_p(\tau_i)(1+\rho_{assay}(i)), \quad i=1,2,...N,
\end{equation}
where the random numbers $\rho_{assay}(i)\in [-b,b]$, $b>0$, are also generated from an almost uniform distribution (Mersenne Twister algorithm). 
The output of the sampling process is the time series defined by the $N$ couples $(t_i, A_i)$, $i=1,2, ..., N$, of times and corresponding measured LH levels. Figure \ref{Construction_A_i} illustrates the construction of the $i$-th sample in a time series.

\begin{figure}[H]
\centering
\includegraphics[width=16.5cm]{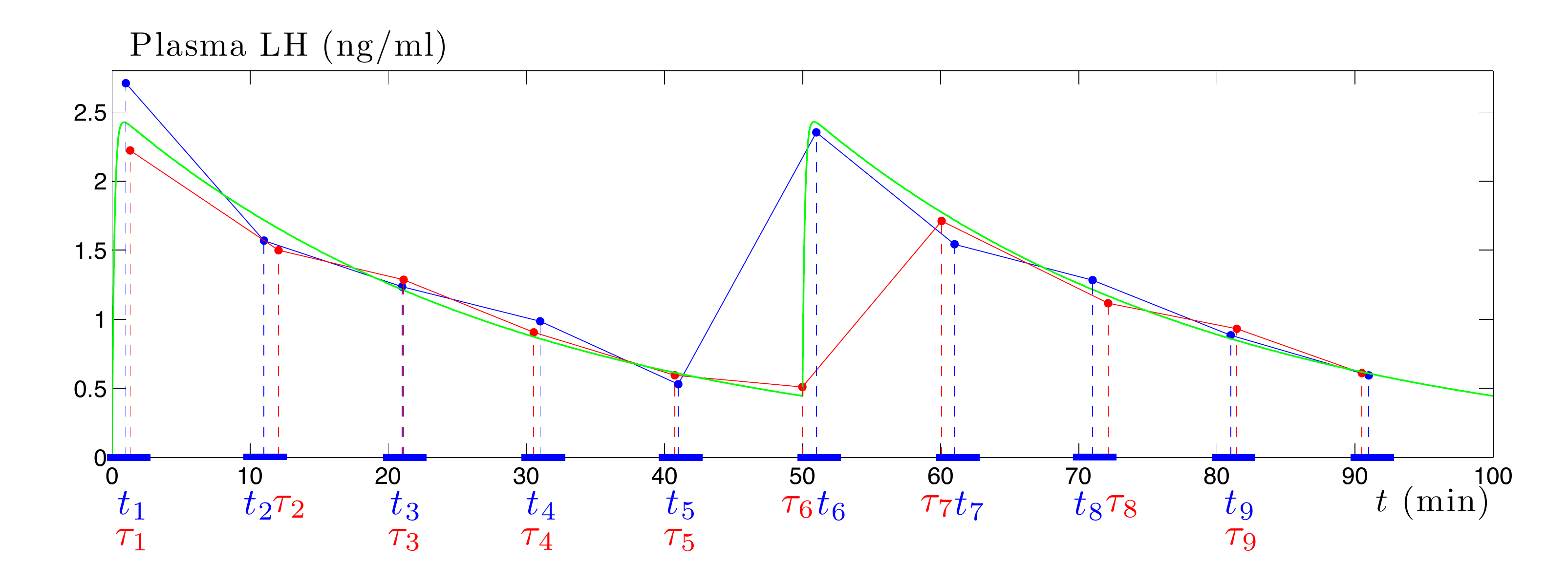}
\caption{
 {\bf  Effect of the variability of the sampling times upon the synthetic time series.}
 The $\tau_i$, $i=1\dots 9$, are the effective sampling times leading to the red-colored LH time series. The $t_i$, $i=1\dots 9$, are the expected sampling times leading to the blue-colored LH time series. The original, non-sampled time series corresponds to the green line. \newline
One can observe an instance of great discrepancy between the LH level measured at time $\tau_6$, which corresponds to the very beginning of the ascending part of a pulse, and the LH level measured at time $t_6$, which corresponds to the maximum of the same pulse.}
\label{t_vs_tau}
\end{figure}

In the context of an experiment, there may be some uncertainty on the exact sampling times (i.e. the precise times at which the samples are retrieved). On the contrary, in our model of the sampling process, we can retrieve the sequence of effective sampling times $(\tau_i)$ for a given synthetic experiment. Figure \ref{t_vs_tau} allows us to visualize the sequence $(\tau_i)$ (red dotted lines) compared to the registered sampling time sequence $(t_i)$ (blue dotted lines). Here, we set the sampling period to $T_s = 10$ min and the shift constant $r=1$ min, hence:
\[
t_i = 10i+1, \quad i=1,2, ..., N
\]
The maximal error is fixed to $15 \%$ of $T_s$, so that $f = 1.5$ min and, for each $i$ from $1$ to $N$, $\tau_i \in [t_i - 1.5, t_i + 1.5]$.

Moreover, Figure \ref{t_vs_tau} compares the results obtained without any variability on the sampling times (blue time series) with those obtained with an error of $\pm 15 \%$ (red time series). It illustrates that this variability can occasionally imply a great difference near a pulse maximum. The 6th blue sample, obtained at $t=t_6=51$ min, corresponds to the theoretical pulse maximum around $2.3$ ng/ml. Yet the 6th red sample, obtained one minute earlier due to the variability of the sampling time, corresponds to the preceding minimal LH level. Consequently, the local maximum of the blue time series, obtained with the 6th sample, is noticeably greater than the local maximum of the red time series, which is obtained with the 7th sample.
\vspace{-0.2cm}

\subsubsection*{Model outputs}

On the endocrinological ground, a LH pulse is an increase in LH blood level triggered by the quick release of LH by the pituitary gland. As illustrated in the preceding section, the moderate clearance rate of LH from the blood underlies the specific asymmetric shape of the pulses, which is characterized by a fast increase immediately followed by a slower decrease. This property has been highlighted in dedicated studies using high frequency sampling (for instance \cite{Keenan_00}: horse, 2 samples per minute) of LH level during a short interval of time.

However, in long-time experiments, the sampling frequency is usually of the order of one per $10$ minutes. Consequently, the precise shape and quantitative properties of the pulses are non longer obvious in the time series. In particular, the theoretical pulse amplitude (theoretical highest level hit during a pulse event) is most of the time not properly reflected by the highest sample obtained during the corresponding event. In the following, we introduce few notions allowing us to differentiate the properties of a theoretical pulse from those of the corresponding pulse obtained from a time series.

The advantages of synthetic time series is that the underlying signal of LH release ($LH(t)$) and the theoretical continuously measured blood LH level ($LH_p(t)$) are available. This corresponds to the ideal experimental situation where one could get high-frequency sampled, variability-free time series retrieved at the same time from the cavernous sinus and jugular blood. With synthetic data, we dispose of reference sets that allow us to identify both LH spikes and pulses without any ambiguity.

Moreover, we can easily test different experimental protocols by changing the value of the parameters  controlling the sampling properties ($T_s ,  r, f, b$ ) and choosing various functions $M_{spike}$ and $P_{spike}$ that determine the time-varying amplitude and frequency of LH spikes released by the pituitary gland.

\begin{figure}[htbp]
\centering
\includegraphics[height=6.9cm]{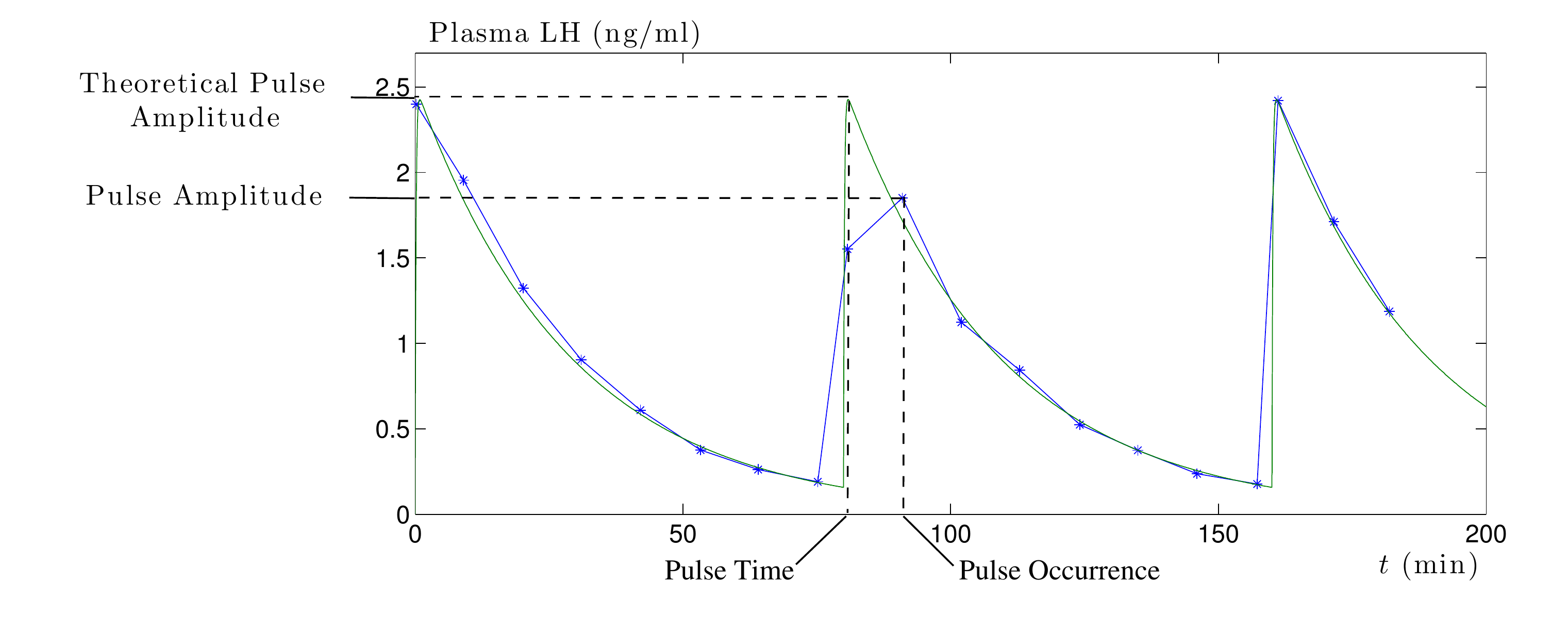}
\caption{
 {\bf Definition of pulse properties in the theoretical case versus experimental case.}
  For a theoretical pulse (i.e. a local maximum in the $LH_p(t)$ signal triggered by a spike in $LH(t)$), we call ``pulse time'' the time at which $LH_p(t)$ admits a local maximum and ``theoretical pulse amplitude'' the value of $LH_p$ at this time. In a time series (either obtained from simulation and synthetic sampling protocol or experimental data), we call ``pulse occurrence'', the time at which the time series admits a local maximum and ``pulse amplitude'' the corresponding value. Both the time values and the amplitude values are different in the theoretical and the experimental cases.}
\label{Vocable}
\end{figure}

\paragraph{Definitions} For sake of clarity, we specify a few notions and terms that will be used in what follows. The definitions are illustrated by Figure \ref{Vocable}. 
For a theoretical pulse (i.e. a peak in the signal $LH_p(t)$ triggered by a spike in $LH(t)$), we define:
\begin{itemize}
\item the theoretical pulse amplitude as the maximal value hit during the event.
\item the theoretical pulse time as the time at which the level hits the theoretical pulse amplitude.
\end{itemize}
For a pulse in a time series (corresponding to a theoretical pulse), we define:
\begin{itemize}
\item the pulse amplitude as the maximal sample obtained during the pulse event,
\item the pulse occurrence as the sample time at which the time series hits the pulse amplitude.
\end{itemize}

\paragraph{Synthetic LH time series obtained from constant spike amplitude and frequency.} We first examine the effects of parameters $r$, $f$ and $b$ in case of a constant spike amplitude $M_{spike}=15$ ng/(ml.min) and constant interspike interval $P_{spike}=100$ min for a same sampling period $T_s=10$ min:
\begin{itemize}
\item Case A: $r=1$ min, $f=0$ min, $b=0$ ;
\item Case B: $r=4$ min, $f=0$ min, $b=0$ ;
\item Case C: $r=4$ min, $f=1.5$ min, $b=0$ ;
\item Case D: $r=4$ min, $f=1.5$ min, $b=10\%$.
\end{itemize}

\begin{figure}[htbp]
\centering
\includegraphics[width=15cm]{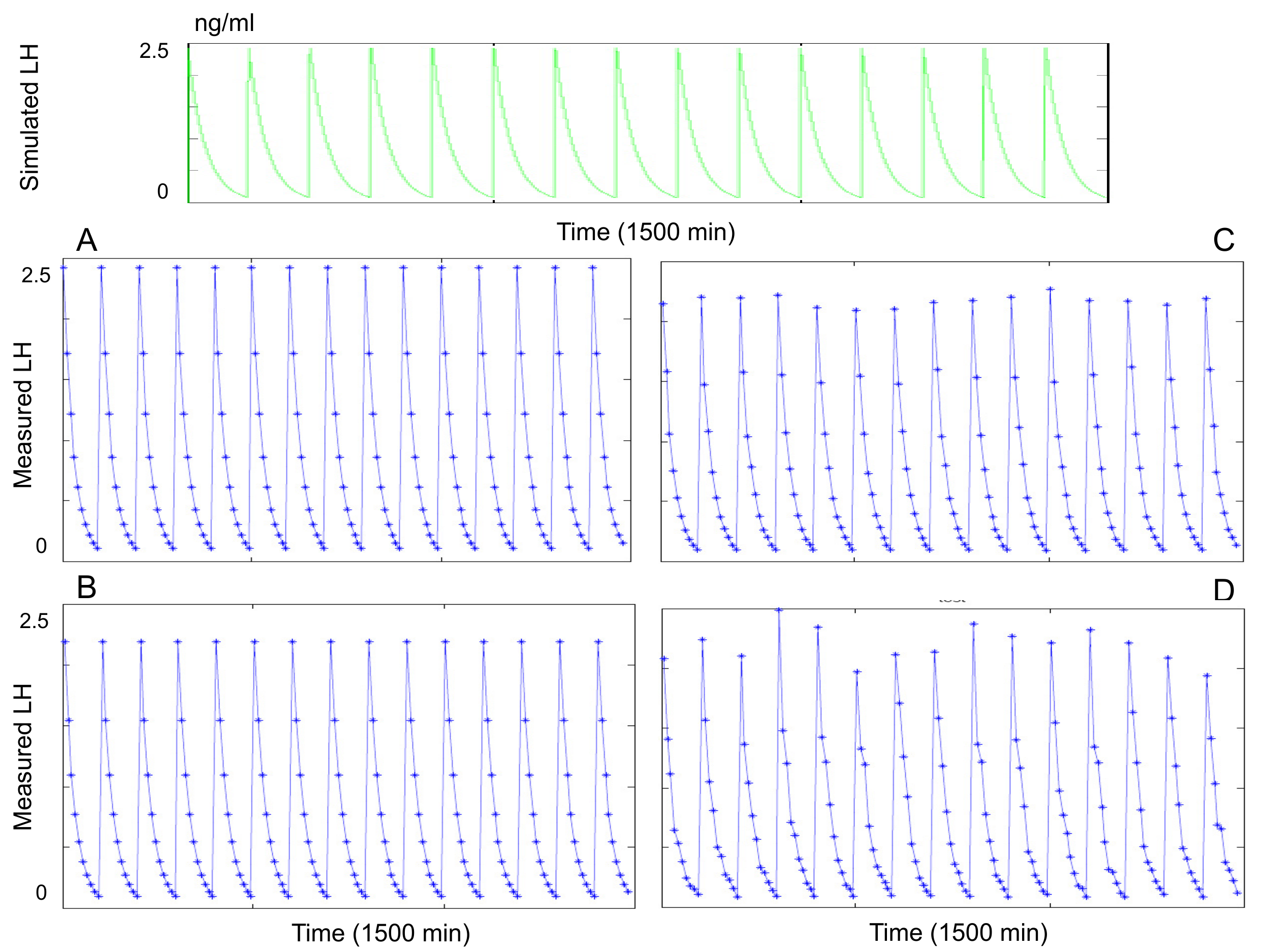}
\caption{
 {\bf  Effect of the sampling process upon a LH level signal with constant amplitude and pulse frequency.}
 In all panels, $M_{spike}=15$ ng/(ml.min), $P_{spike}=100$ min, $T_s=10$ min.
Top panel: theoretical continuously measured LH blood level. Panels A, B, C, D: sampling points (blue stars) of the time series obtained from the top panel signal through the sampling protocol. Panel A: first sampling time at $r=1$ min, without any variability in the sampling process. Panel B:  first sampling time at $r=4$ min, without any variability in the sampling process. Panel C: first sampling time at $r=4$ min, with variability in the sampling times ($\pm 1.5$ min). Panel D: first sampling time at $r=4$ min, with variability both in the sampling times ($\pm 1.5$ min) and the assays ($\pm 10 \% $).}
\label{Sig_Freq_Cst}
\end{figure}

The top panel of Figure \ref{Sig_Freq_Cst} displays the solution $LH_p(t)$ (green curve) of equation \eref{LHplasma}, i.e. the theoretical continuously measured LH blood level. Each panel from A to D shows the time series (blue stars) obtained through the sampling protocol for the values of parameters $r$, $f$ and $b$ specified above.

\begin{figure}[H]
\centering
\includegraphics[width=15cm]{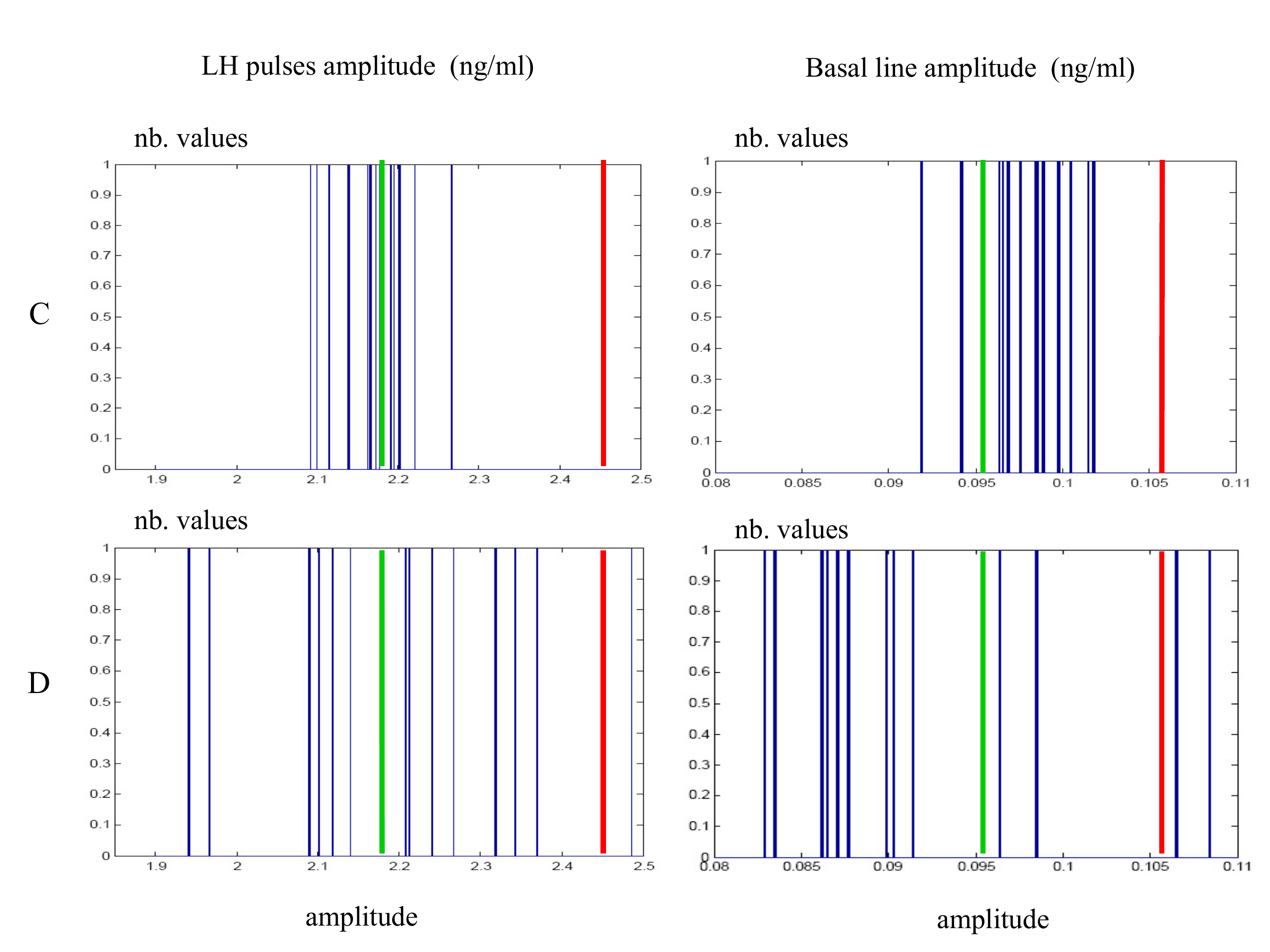}
\caption{
{\bf Histograms of LH pulses and basal line amplitudes in simulated LH series with constant amplitude and pulse frequency.}
Left panels correspond to the maximal value (pulse amplitudes) distribution whereas right panels correspond to the minimal value (basal line amplitude) distribution measured from the four cases of figure \ref{Sig_Freq_Cst}. The A and B time series, that only differ in the first sampling time, display constant (yet different) pulse amplitude. Red bars stand for case A ($r = 1$ min) value of the pulse (2.425) and basal line (0.107) amplitudes. Green bars stand for case B ($r=4$ min) value of the pulse ($2.188$ ng/ml) and basal line ($0.096$ ng/ml) amplitudes. Blue bars stand for case C ($r=4$ min; $f=1.5$ min) and case D ($r=4$ min; $f=1.5$ min; $b= 10 \%$, i.e. a variability of $\pm 10 \%$ in the LH assays) distributions of LH pulse and basal line amplitudes. 
Case D pulse ($\max=2.486$ ng/ml; $\min=1.940$ ng/ml) and basal line ($\max=0.108$ ng/ml; $\min=0.082$ ng/ml) amplitude distributions are wider than case C pulse ($\max=2.266$ ng/ml; $\min=2.092$ ng/ml) and basal line ($\max=0.101$ ng/ml; $\min= 0.092$ ng/ml)  amplitude distributions, due to combined variabilities in the sampling times and assays.}
\label{Hist_Freq_Cst}
\end{figure}

Figure \ref{Hist_Freq_Cst} details, for each time series obtained in case C or D, the distribution of LH pulse amplitude (blue bars in left panels), i.e. local maxima, and LH levels at the basal line (blue bars in right panels), i.e. local minima. For sake of comparison, the constant amplitudes obtained in cases A and B have been marked with red and green bars respectively.

In cases A and B ($f=b=0$), one obtains a strictly periodic pattern of sampled LH levels since $P_{spike}$ is a multiple of $T_s$. However, depending on the beginning of the sampling process  $r$ ($1$ min in time series A and $4$ min in time series B), the maximum sample value varies from $2.425$ ng/ml in case A (red bars in left panels of Figure \ref{Hist_Freq_Cst}) to $2.188$ ng/ml in case B (green bars). This shows the importance of the phase between the pulsatile LH signal and the periodic sampling process in the resulting observed pulse amplitude. It is worth noticing that this phase cannot be controlled at all in experimental conditions since the delay elapsed from the last pulse time is not known at the beginning of the experimental sampling process. The dependence of the basal line on this phase is weaker but still exists: it varies from $0.107$ ng/ml in time series A to $0.096$ ng/ml in time series B.

By comparing panels B and C of Figure \ref{Sig_Freq_Cst}, we can observe the impact of the variability in the sampling times ($f=0$ min and $f=1.5$ min in time series B and C respectively). With variable sampling times, the pulse amplitude along time series is also variable, although the original continuous signal $LH_p(t)$ is perfectly periodic. This variability is shown in the top panels of Figure \ref{Hist_Freq_Cst}: the various LH pulse (resp. basal line) amplitudes obtained in case C (blue bars) are scattered around the constant case B pulse (resp. basal line) amplitudes (green bar).

The impact of the variability in the assays is illustrated by the enhanced dispersal of the LH pulse and basal line amplitudes obtained in case D (blue bars in bottom panels of Figure \ref{Hist_Freq_Cst}) for which $b=10 \%$ compared to the case C amplitudes ($b=0$).

In any case of either shifted (case B) or noised (cases C and D) time series, the pulse amplitudes are undervalued with respect to the genuine amplitude (correctly assessed only in case A). It may nevertheless happen that the effective sampling time coincides with a (genuine) pulse time. In that case, the pulse amplitude can be overvalued if the sign of the assay error is positive (instance of the blue bar on the right of the red bar in the left panel of case D).

\paragraph{Synthetic LH time series obtained from time-varying spike amplitude and frequency.} We now further examine the effects of the sampling process upon theoretical continuously measured LH level with time-varying pulse amplitude and frequency:
\begin{itemize}
\vspace{-0.1cm}
\item Case E: constant spike amplitude $M_{spike}=15$ ng/(ml.min) and decreasing interspike interval $P_{spike}(t)=100-\frac{t}{30}$ (expressed in min);
\vspace{-0.1cm}
\item Case F: decreasing spike amplitude $M_{spike}(t)=15-8.7 \ 10^{-3}t$ (expressed in ng/(ml.min)) and decreasing interspike interval $P_{spike}(t)=100-\frac{t}{30}$ (expressed in min).
\end{itemize}
\vspace{-0.1cm}
The sampling period is the same in both cases: $T_s=10$ min.

\begin{figure}[htbp]
\centering
\includegraphics[width=16cm]{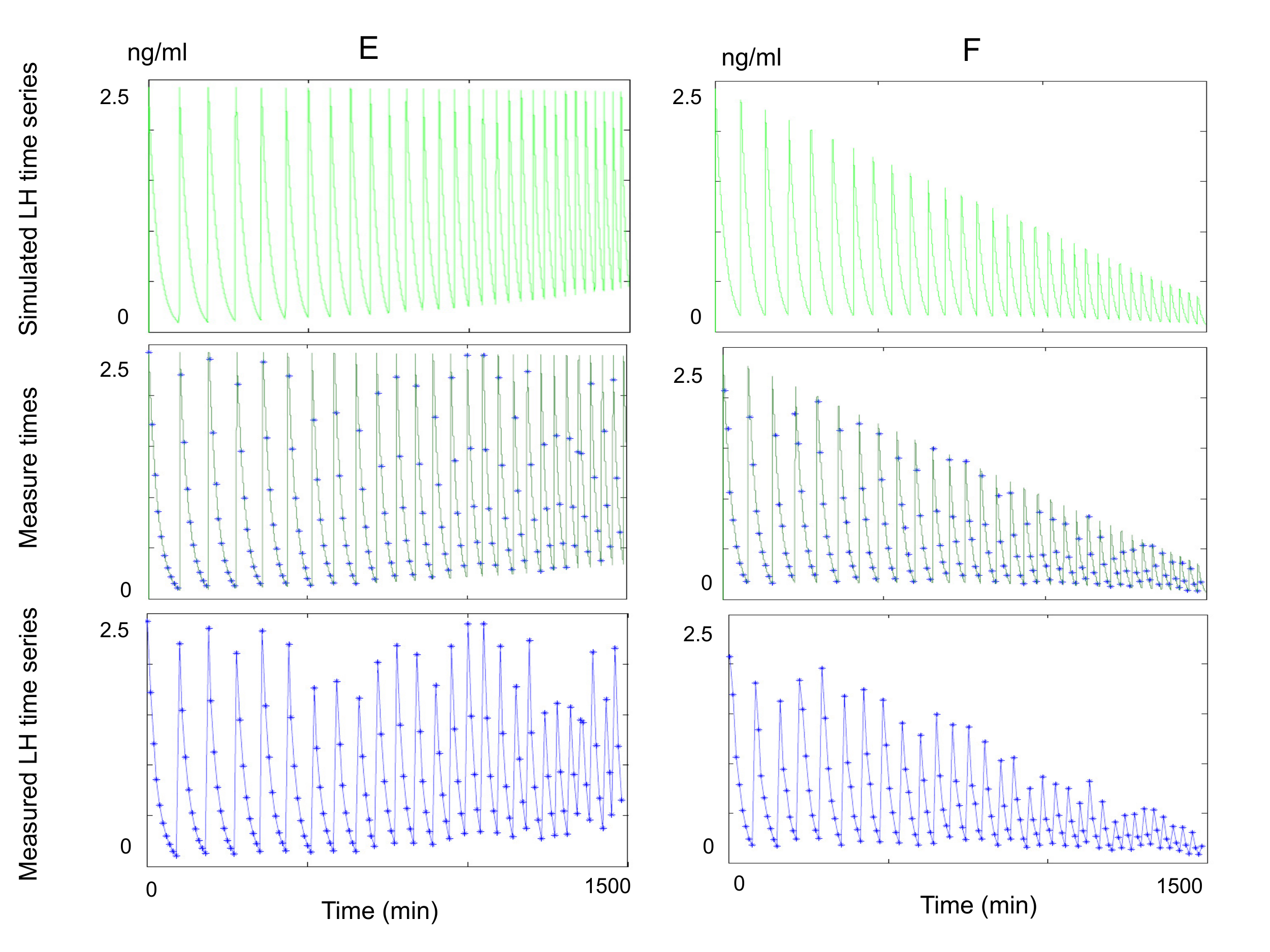}
\caption{
{\bf  Effect of the sampling process upon a LH level signal with regular increasing sampling frequency.}
Case E: the pulse amplitude remains almost constant and the basal line regularly increases. Case F: the pulse amplitude regularly decreases and the basal line regularly decreases. The top panels represent the fine step simulation of plasma LH level. The middle panels represent the time series (blue stars) along the theoretical continuously measured plasma LH level (green curve). The bottom panels represent the resulting LH measured time series (measured LH levels versus sampling times linked with segments). In both cases, the sampling period is $T_s = 10$ min. In case E, the initial time sampling occurs at the first minute of the simulation ($r=1$ min), without any sampling time ($f=0$ min) or assay ($b=0$) variability. In case F, the initial time sampling occurs at the fourth minute of the simulation (r=4 min), with sampling time ($f=2$ min) and assay ($b=5 \%$) variability.}
\label{Sig_Freq_Incr}
\end{figure}

\begin{figure}[htbp]
\centering
\includegraphics[width=16cm]{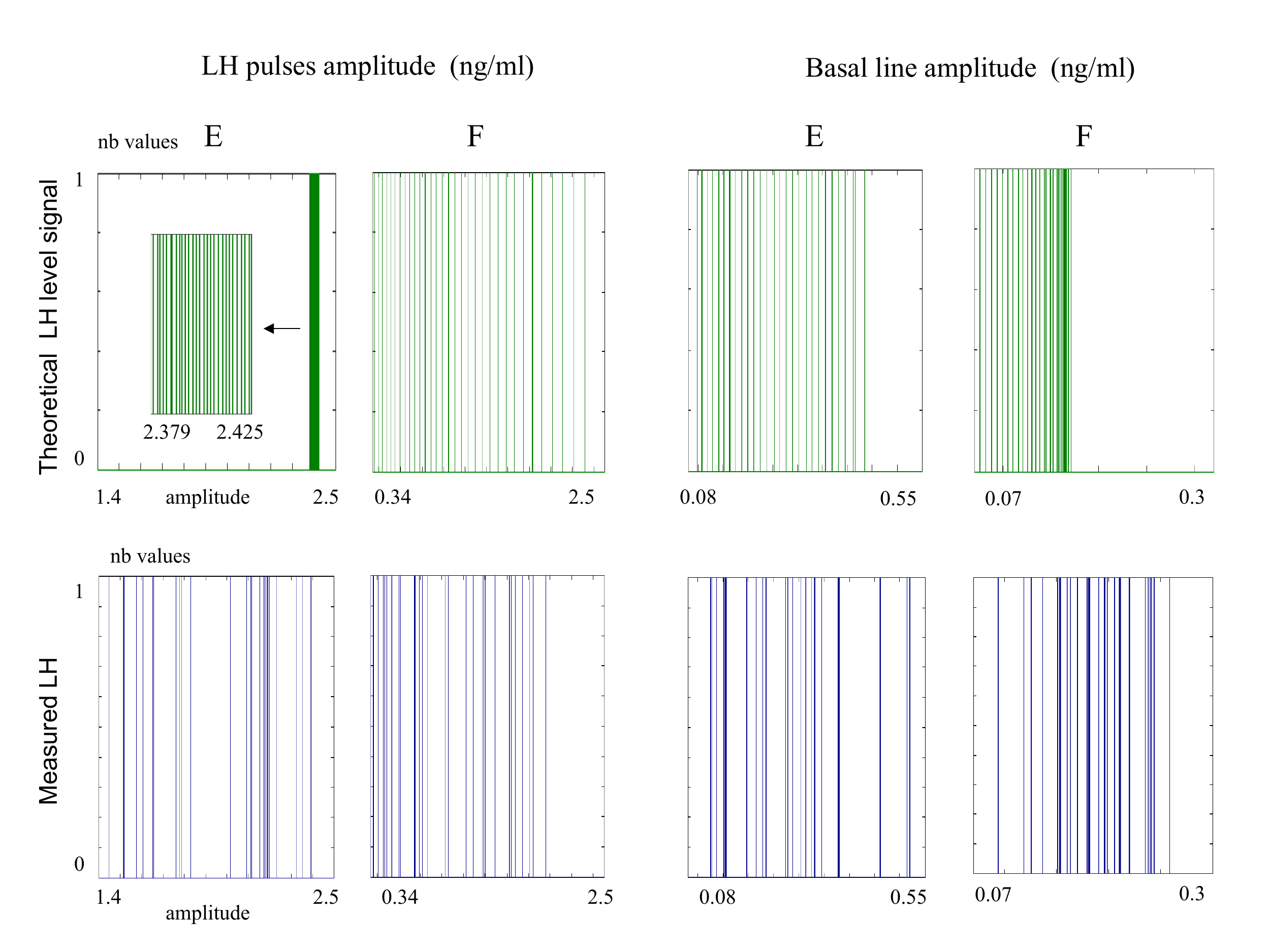}
\caption{
{\bf  Histograms of LH pulses and basal line amplitudes in simulated LH series with increasing pulse frequency.}
Left panels correspond to the maximal value (pulse amplitude) distributions whereas right panels correspond to the minimal value (basal line amplitude) distributions, measured from the two cases of figure \ref{Sig_Freq_Incr}. Top panels (green lines) display the initial LH plasma level distributions whereas the bottom panels (blue lines) display the resulting LH time series. A zoom on the amplitude distribution is shown as an insert in case E of the top left panel (left arrow).
While the distributions are regular in the theoretical time series, they become completely irregular in the sampled time series. As a result, the range of amplitudes is shortened.
Regarding the pulse amplitude distribution in case F, it is worth noticing that all measured values ($\max=1.959$ ng/ml, $\min=0.302$ ng/ml) are lower  than the corresponding theoretical values ($\max=2.315$ ng/ml, $\min=0.353$ ng/ml). In case E, the theoretical pulse amplitudes vary from $2.425$, to $2.379$ ng/ml whereas measured pulse amplitudes vary from $2.395$ to $1.447$ ng/ml.
On the contrary, regarding the basal line amplitude distribution, it is worth noticing that all the measured values (E: $\max=0.519$ ng/ml, $\min=0.125$ ng/ml; F: $\max=0.258$ ng/ml, $\min=0.094$ ng/ml) are greater than the theoretical minimal values (E: $\max= 0.434$ ng/ml, $\min= 0.098$ ng/ml; F: $\max= 0.171$ ng/ml, $\min= 0.075$ ng/ml).}
\label{Hist_Freq_Incr}
\end{figure}

Figure \ref{Sig_Freq_Incr} gives some examples of time series obtained in cases E and F. In each case, the green curve represents the solution $LH_p(t)$ (theoretical continuously measured plasma LH level) and the blue stars represent the time series obtained through the synthetic sampling protocol.
Figure \ref{Hist_Freq_Incr} details, for cases E and F, the distribution of LH pulse amplitudes (left panels) and successive LH levels at the basal line (right panels) both for the theoretical continuously measured LH blood level (top panels, green bars) and the time series obtained through the sampling protocol (bottom panels, blue bars).

In case E, the spike frequency increases from $1$ spike per $100$ min to $1$ spike per $50$ min. Hence, the spike release arises more and more often along time, so that the LH blood pulses are successively triggered from a higher and higher basal level. Consequently, the basal line and, in a lesser extent, the pulse amplitude of the theoretical LH level undergo a small and smooth increase. Moreover, the number of samples per pulse decreases drastically as the pulse frequency increases, so that the LH time series looks noisier at the end than at the beginning of the time series (Figure \ref{Sig_Freq_Incr}). This effect implies that the pulse amplitudes are spread out by the sampling protocol much stronger in case E (from $2.395$ to $1.447$ ng/ml) than in cases C and D (Figure \ref{Hist_Freq_Incr}). Hence the time variations in the pulse frequency enhances the variability brought about by the sampling process.

Case F represents a situation that is naturally encountered in the physiological dynamics of LH secretion: the same increase in the spike frequency as in case E, combined with a decrease in the amplitude. As in the preceding cases, each of the LH pulse amplitude in the time series is smaller than the corresponding one in the theoretical LH level. Additionally in case F, the amplitude at the basal line is sensibly raised up by the sampling protocol. Hence, the difference between the pulse amplitude and the basal level is strongly deprecated, which adds to the noisy character of the ending of the time series.

\vspace{-0.15cm}
\subsection*{Pulse detection algorithm}

\vspace{-0.1cm}
In the context of automatic pulse detection in a time series, a pulse is a peak (i.e. a local maximum of the time series) fulfilling given criteria. The most challenging issue in the algorithm design consists in formalizing the biologically relevant criteria that discriminate the pulses from other peaks. Among these criteria, the amplitude is the most obvious. However, as illustrated in the preceding section, it cannot be used as an infaillible criterium for automatic pulse detection since the quantitative features of a pulse (absolute amplitude, amplitude from baseline, ...) are really altered in an unknown way by the sampling process.

For sake of robustness and acuteness of our pulse detection algorithm, we have introduced a selection process based on multiscale criteria involving different properties of the peaks. Besides the series of pulse occurrences, the main output of our algorithm is the series of the corresponding InterPulse Interval (IPI) together with a tunnel of confidence (IPI tunnel) related to the regularity of the pulse frequency variations.

\subsubsection*{Algorithm foundations} \label{AlgoFound}

\paragraph{Notations and definitions}
We consider a time series of $N$ points obtained with a $T_s$-periodic sampling process and we note the sampling times:
\[
t_k=(k-1)T_s, \quad k= 1, 2, \dots N.
\]
Hence, the $k^{th}$ sampling time is $t_k$ (in particular, the first sampling time is $t_1=0$). We note either $A_k$ or $A(t_k)$ the value corresponding to the $k^{th}$ sample and call $k$ a ``time index''. For sake of simplicity in the following explanations, we will refer either to the time index $k$ or the corresponding time $t_k$.

The pulse detection algorithm is based on a sequence of different processes. Some of them consist in researching high amplitude peaks that can potentially be classified as pulses, others aim to remove, among the formerly selected peaks, those that do not fit other properties met by genuine pulses. In the following, we note $P=(p_i)_{1 \leq i \leq s(P)}$ the vector storing dynamically the time indexes at which the algorithm detects the summit of a potential pulse. $s(P)$ is the size (number of elements) of $P$. The algorithm modifies vector $P$ in such a way that indexes $p_i$ are always sorted in increasing order. Hence, at any time along the algorithmic process, $s(P)$ pulses are detected, $t_{p_i}$ is the occurrence of the $i^{th}$ detected pulse and $A_{p_i}$ is the amplitude of the $i^{th}$ detected pulse.

\begin{figure}[!ht]
\centering
\includegraphics[height=6.7cm]{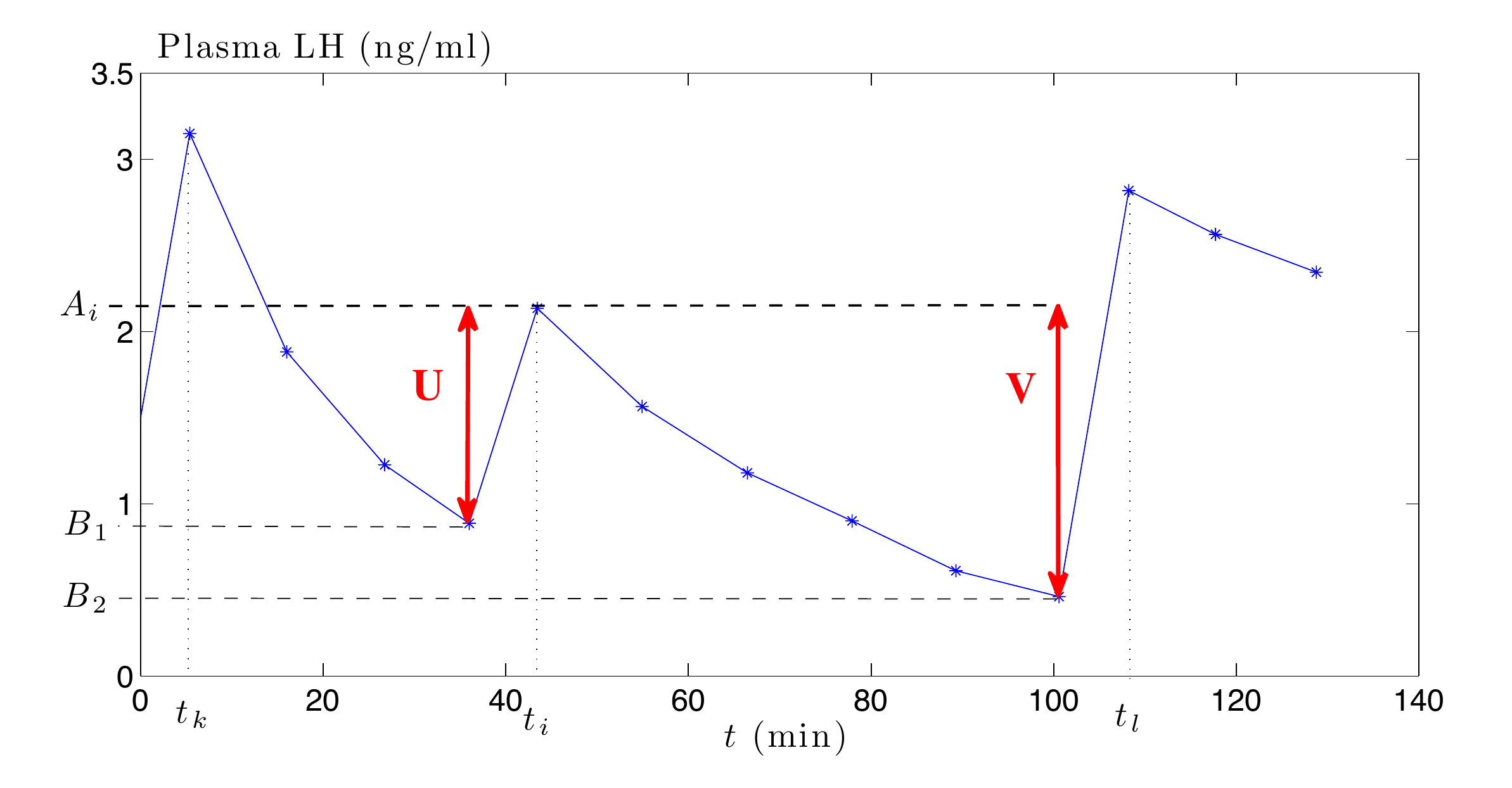}
\caption{
 {\bf Definition of the magnitude of a peak.} Let $P$ be a given vector of potential pulses. Considering the peak occurring at $t_i=43$ min, we assume that $P$ contains the pulses just before and after $t_i$ occurring at $t_k=4$ min and $t_l=108$ min. We define $U$ (resp. $V$) as the difference between the peak amplitude ($A_i=2.2$ ng/ml) and the minimum value $B_1$ (resp. $B_2$) of the time series between $t_k$ and $t_i$ (resp. $t_i$ and $t_l$) : $B_1=0.8$ ng/ml (resp. $B_2=0.4$ ng/ml). The magnitude of the peak occurring at $t_i$ is the geometric mean between $U$ ($1.4$ ng/ml) and $V$ ($1.8$ ng/ml). Here, the peak magnitude is equal to $1.587$ ng/ml.}
\label{Magnitude}
\end{figure}

In order both to moderate the importance of the amplitude and to account for several characteristics of the pulse shape, we need to introduce the notions of ``height'' and ``magnitude'' of a peak as well as that of ``relative magnitude'' between two peaks. Let us suppose that the $i$-th sample of the time series (occurring at time $t_i$) corresponds to a peak (i.e. $A_{i-1}<A_{i}$ and $A_{i+1}<A_{i}$). The height of this peak is defined as the difference between its amplitude $A_{i}$ and the lowest value of $A_j$ within the sampled time series, denoted by $\underline{A}$, i.e.:
\[
H_i=A_{i}-\min_{j=1..N} A_j=A_i-\underline{A}.
\]
Let us assume additionally that a set of potential pulses $P$ is identified from the time-series and the closest pulses registered in $P$ before and after $t_i$ occur at $t=t_k$ and $t=t_l$ respectively (see Figure \ref{Magnitude}). We define the two minimal values of the time series between $t_k$ and $t_i$ on one hand, and between $t_i$ and $t_l$ on the other hand:
\begin{align*}
     B_1 &= \min_{j\in \{ (k+1),(k+2),\dots, (i-1)\}}  A_k \\
     B_2 &= \min_{j\in \{ (i+1),(i+2),\dots, (l-1)\}}  A_k
\end{align*}
Then, we define the magnitude of the peak corresponding to the $i$-th sample as the geometric mean of $A_{i}-B_1$ and $A_{i}-B_2$, i.e. $\sqrt{(A_i-B_1)(A_i-B_2)}$.

Finally, let us consider two peaks occurring at $t_k$ and $t_l$ with $t_k<t_l$. Let us call $B_0$ the minimal value of the time series between $t_k$ and $t_l$:
\[
B_0 = \min_{j\in \{ (k+1),(k+2),\dots, (l-1)\}}  A_k
\]
We call ``relative magnitude'' between the two peaks the geometric mean between $A_{k}-B_0$ and $A_{l}-B_0$.

The notion of peak height does not depend on the vector of potential pulses. It represents a normalization of the amplitude among the time-series with respect to the lowest value of the time serie. The peak magnitude can change as the vector $P$ of identified pulses evolves and, consequently, will be used to compare a peak with its direct neighbors. The interest of the magnitude is to take into account the local baseline, without being sensitive to differences in the baseline from one side or the other of a peak. The relative magnitude between two peaks gives a semi-local reference for the pulse magnitude. In particular, the magnitude of a pulse can be usefully compared to the relative magnitude between its direct neighbors (i.e. the pulses just before and after it).

\paragraph{Pulse selection process}
\begin{itemize}
\item \textbf{Initialization:} We first fill vector $P$ by selecting time indexes corresponding to great height samples. Even if, at this stage, we intend to recover a large enough set of potential pulse indexes, the time intervals between pulses should be consistent with the maximal frequency. Accordingly, we introduce a parameter $T_{p}$, called the nominal period, defined as the smallest time duration in which, from one pulse occurrence, one expects the following one. \\
Starting from the time index $p_1$ of the maximal sample in $[0, 2T_p]$, we locate the time index $m_2$ of the minimal sample in $[t_{p_1},t_{p_1}+T_p]$ and then we retrieve the time index $p_2$ of the maximal sample in $[t_{m_2}, t_{m_2}+T_p]$. We iterate the process along the whole time series to obtain the initial guess of potential pulse indexes $P=(p_i)$. \\
This method is a trade-off between selecting all the peak indexes in the time-series (with the drawback that there will be too many of them if the time series is noisy) and selecting only local maxima corresponding to great amplitudes (with the drawback that there will be too few of them if the time series is smooth and the pulse frequency is low).
\item \textbf{Remove too small peaks:} Once vector $P$ is initialized, some of the registered indexes may correspond to small sample values. As illustrated in the first section, the pulse amplitude is strongly altered by the sampling process and, consequently, it cannot be used as an infallible criterium to select the pulses. Hence, we use multi-scale criteria based on the notions of height and magnitude to determine which indexes corresponds to too small peaks and to remove them from vector $P$.
\begin{itemize}
\item \textbf{Global relative criterium:} We aim to remove peaks whose height is small in comparison with the height of all detected pulses. We define the ``median height'' as the median of the detected pulse heights. For each pulse, if the ratio between its own height and the median height is less than a threshold parameter $\lambda_r$, it is removed from $P$. \\
We have chosen the median instead of the (arithmetic or geometric) mean since it is less sensitive to the presence of great height peaks..
\item \textbf{Semi-local relative criterium:} We compare the magnitude of each peak with the relative magnitude between the immediately preceding and following pulses. If this comparison is not conclusive with respect to the threshold $\lambda_r$ introduced above, we remove the peak index from vector $P$. \\
It is worth noticing that the geometric mean used to define the magnitudes provides robustness (compared to arithmetic mean, for instance) to this criterium with respect to possible local variations of the base line.
\item \textbf{Global absolute criterium:} We compare the magnitude of each pulse to an absolute threshold $\lambda_a$. If the comparison is not conclusive, the corresponding time index is removed from vector $P$. \\
This criterium precludes non significant elevations in the baseline to be considered as potential pulses. Hence, parameter $\lambda_a$ corresponds to the assay detection threshold.
\end{itemize}
At this level of the selection process, vector $P$ only contains the time indexes corresponding to peaks with sufficiently great height and magnitude (with respect to threshold $\lambda_r$ and $\lambda_a$) to be classified as pulses. However, as the initial guess for $P$ may have skipped some potential pulses, the next step of the selection process consist in retrieving the missed pulses.
\item \textbf{Retrieve missed pulses:} Between each pair of successive pulses registered in $P$, we examine each peak. If this peak fulfills the semi-local relative criterium, the corresponding time index is added to vector $P$. \\
By construction, such a retrieved peak almost automatically fulfills the other two global criteria.
\item \textbf{Shape-based criterium:} The pulse duration has to be consistent with available knowledge on pulse half-life. For a fixed sampling rate, a detected pulse should extend over a minimum number of consecutive experimental data. Consequently, we intend to remove what we call ``3-point peaks'' for which the immediately preceding and following samples are local minima of the time series (see top panel of Figure \ref{3-point-peak}).
\vspace{-1.5cm}
\begin{figure}[H]
\centering
\includegraphics[width=13cm]{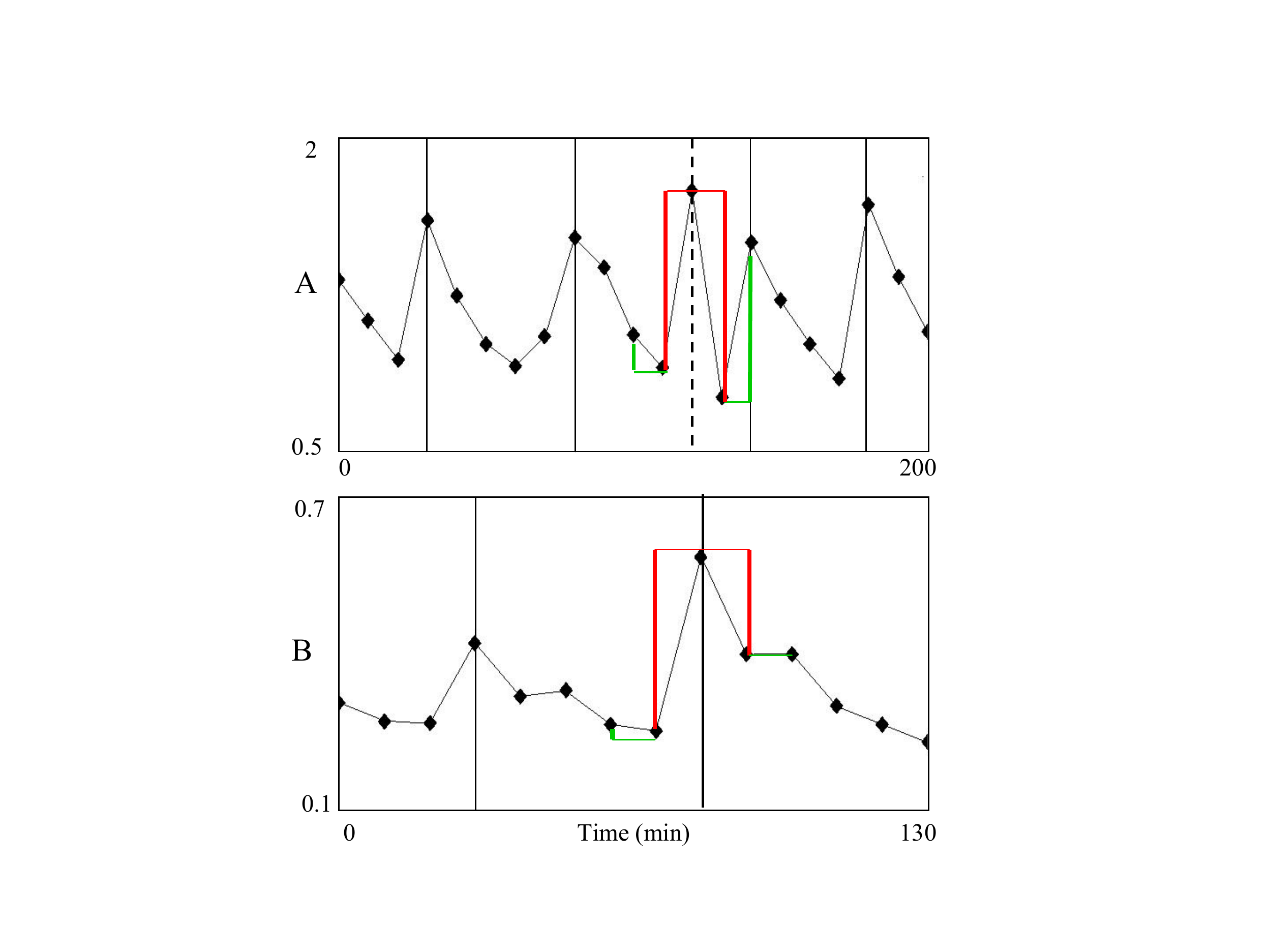}
\vspace{-1cm}
\caption{
{\bf  Identification of a 3-point peak pattern.}
Parameter $\lambda_{3p}$ corresponds to the ratio between the arithmetic mean of the amplitude of the neighboring points of rank 2 (green lines) and the geometric mean of the amplitude of the immediate neighbors (red lines). Panel A: the selected pulse (dashed vertical line) is identified as a genuine 3-point peak. Panel B: the selected pulse (solid vertical line) is not identified as a 3-point peak, since it belongs to a genuine, asymmetric LH pulse with an exponential decrease, albeit locally noised.
}
\label{3-point-peak}
\end{figure}
\vspace{-0.5cm}
However, due to possible noise in the time-series, a pulse may appear as a 3-point peak. But, in this case, the pulse is expected to be not ``too sharp'' (see bottom panel of Figure \ref{3-point-peak}). Thus we remove from $P$ the time indexes corresponding to peaks with a ``sharpness coefficient'' greater than a chosen threshold $\lambda_{3p}$. \\
The precise definition of ``sharpness coefficient'' will be detailed in Step 5 of the algorithm (see next subsection ``Algorithmic pulse detection''). We only enlighten here that this criterium is based on the asymmetric shape of a pulse and allows one to get rid of the genuine peaks produced by occasional experimental errors.
\end{itemize}

\paragraph{Vector of InterPulse Interval (IPI) and IPI tunnel} \label{IPIoutliers}

At a given step, we define the vector of InterPulse Intervals (IPI) from the current $P$ vector by:
\[
\Theta = (\theta _i)_{1\leq i \leq s(P)-1} \text{ with }\theta _i=t_{p_{i+1}}-t_{p_{i}}.
\]
As we aim to apply the algorithm mainly to hormonal time series, we designed a process to take into account some degree of regularity in the rate of change in the pulse frequency. Hence, given a vector $P=(p_i)$ of identified potential pulses, we introduce a cubic function $\phi (i)$ fitted to $\theta _i$ in a least squares sense. Function $\phi $ gives an averaged, yet time evolving representation of the IPI built from the global sequence $\Theta $.

Then, we build a tunnel in $[t_{p_2}, t_{p_{s(P)}}] \times \mathbb{R_+}$ delimited by the graphs of two piecewise linear functions of time $\psi_{\inf}$ and $\psi_{\sup} $ built from function $\phi $. Precisely, for each pulse time $t_{p_i}$ ($i=2...s(P)$), we define $\psi_{inf}(t_{p_i})=(1-\alpha) \phi(i-1)$ and $\psi_{sup}(t_{p_i})=(1+\beta) \phi(i-1)$ to draw the lower and upper boundaries of the tunnel. Thus, parameters $\alpha $ and $\beta $ tune the width of the so-called ``IPI tunnel'' ; they represent a quantification of the pulse frequency regularity. The tunnel allows one to assess the regularity of the IPI time variations and can help the user to (i) classify a specific pulse as a potential outlier with respect to the pulse frequency properties of the series, (ii) identify and localize a possible rupture in the secretion rhythm.

To illustrate the use of the IPI tunnel, let us consider a time series for which the sequence of pulse indexes $Q=(q_i)$ is easy to find by sight. We assume that the time series displays regular pulsatility, i.e. the pulse frequency undergoes smooth variations along time. Under this condition, function $\phi $ is a good approximation of the IPI sequence.

Let us first consider the case of a lack of detection: let $P$ be the vector formed by the pulse indexes in $Q$ except one (the $i_0^{th}$). In the course of the automatic pulse detection, this case may happen if the maximum amplitude corresponding to this pulse is low. Then, the corresponding IPI sequence $(\theta ^P_i)$ is given by:
\begin{equation}
\theta ^P_i =
\begin{cases}
\theta ^Q_i \quad \text{ if } 2 \leq i \leq i_0-1, \\
\theta ^Q_{i_0} + \theta ^Q_{i_0+1} \quad \text{ if } i=i_0, \\
\theta ^Q_{i+1} \quad \text{ if } i_0+1 \leq i \leq s(Q)-1=s(P).
\end{cases}
\end{equation}
Under the regularity assumption, each IPI $\theta _i$ should be close to (or even in) the range delimited by the values of its neighbors $\theta _{i-1}$ and $\theta _{i+1}$. On the contrary, in the case of vector $P$, the $(i_0-1)^{th}$ IPI ($\theta ^P_{i_0}$) is noticeably greater than the maximum of its neighbors. More precisely, its value is twice the mean of its neighbors and, consequently, is approximately twice the $\phi $ value.

Now, let us consider the case of an over-detection: let $P$ be the vector formed by the pulse indexes in $Q$ plus an extra pulse occurring at time $t=t_k$ lying between the $(i _0-1)^{th}$ and the $i _0^{th}$ pulse times stored in vector $Q$. Hence, $q_{i_0-1}<k<q_{i_0}$. Then: 
\begin{equation}
\theta ^P_i =
\begin{cases}
\theta ^Q_i \quad \text{ if } 2 \leq i \leq i_0-1, \\
\theta ^Q_{i-1} \quad \text{ if } i_0+2 \leq i \leq s(Q)+1=s(P),
\end{cases}
\end{equation}
and moreover:
\[
\theta ^P_{i_0} + \theta ^P_{i_0+1}=\theta ^Q_{i_0}.
\]
Under the regularity assumption, either $\theta ^P_{i_0}$ or $\theta ^P_{i_0+1}$ is too small compared to the expected IPI range  delimited by $\theta ^P_{i_0-1}=\theta ^Q_{i_0-1} $ and  $\theta ^P_{i_0+2}=\theta ^Q_{i_0+1}$. In the less discriminating case,  the extra pulse lies close to the middle of its neighbors $(q_{i_0-1}+q_{i_0})/2$. Then, $\theta ^P_{i_0}$ and $\theta ^P_{i_0+1}$ are almost equal to the half of the expected IPI $\theta ^Q_{i_0}$. Hence, even if several combinations of $\theta ^P_{i_0}$ and $\theta ^P_{i_0+1}$ values exist, depending on the position of time $k$ compared to $q_{i_0-1}$ and $q_{i_0}$, one of the IPIs is always less than the half of the $\phi $ value, which indicates a potential over-detection.

Finally, an appropriate choice of the values of $\alpha $ and  $\beta $ allows one to discriminate the pulses that break the frequency regularity and indicate a rupture in the secretion rhythm.

\subsubsection*{Algorithmic description} \label{Algo}

\paragraph{Overview of the algorithm.}
Besides the time series of $N$ samples and the sampling period $T_s$, the user is asked to enter values of four parameters that rule the acuteness of the detection process and two parameters that define the width of the IPI tunnel. We use the parameter names in the following algorithm description. We give a tutorial in the result section for helping the users to choose the appropriate values.
\begin{enumerate}
\item The nominal period $T_p$ represents the smallest duration in which, from one pulse occurrence, one expects the following one. In what follows, we note:
\begin{equation} \label{kp}
k_p = \left\lfloor \frac{T_p}{T_s} \right\rfloor.
\end{equation}
Note that the time interval $[t_{k+1}, t_k+T_p]$ contains the $k_p$ time indexes $\{k+1, ..., k+k_p\}$.
\item The relative magnitude threshold: $\lambda_{r}$,
\item The absolute magnitude threshold ratio: $\lambda_{a}$,
\item The relative magnitude threshold ratio for 3-point peaks: $\lambda_{3p}$
\item The ratios defining the edges of the IPI tunnel: $\alpha $, $\beta $.
\end{enumerate}

The algorithm body is segmented into 6 steps. The third step is itself segmented into 3 substeps. Each step is described in detail below with an accompanying box enunciating the corresponding pseudo-code.

The output of the algorithm consists in:
\begin{enumerate}
\item The set of detected pulse occurrences ($P$),
\item the corresponding sequence of IPIs ($\Theta $),
\item the lower and upper edges of the IPI tunnel (respectively $\psi_{\inf}$ and $\psi_{\sup}$),
\item the sets of lower and upper outliers ($P_{low}$ and $P_{up}$).
\end{enumerate}
and associated graphical outputs (of which we make an intensive use in the result section):
\begin{enumerate}
\item the sampled time series $(A_i)_{i=1}^n$ versus the sampling times with identified pulse occurrences,
\item the graph of the IPI versus the corresponding pulse occurrences (more precisely: $(p_{i+1}, \theta_i)_{i=1}^{s(P)-1}$) together with the IPI tunnel.
\end{enumerate}

\begin{figure}[H]
\fbox{\begin{minipage}{\textwidth}
\begin{algorithmic}
\STATE \textbf{Input data:} $(A_i)_{i=1}^n$, $T_s$, $T_p$, $\lambda_r$, $\lambda_a$, $\lambda_{3p}$, $\alpha$, $\beta$.
\STATE \hrulefill
\STATE Algorithm step 1
\STATE ~~~~~~~~~$\vdots$
\STATE Algorithm step 6
\STATE \hrulefill
\STATE \textbf{Output data:} $P$, $\Theta $, $\psi_{\inf}$, $\psi_{\sup}$, $P_{low}$, $P_{up}$.
\vspace{0.3cm}
\STATE \textbf{Related graphical outputs:} 
\begin{itemize}
\item the time series $(t_i,A_i)_{i=1}^n$ with vertical bars indicating the detected pulse occurrences,
\item the graph of the IPI versus the corresponding pulse occurrences and the IPI tunnel: $(p_{i+1}, \theta_i)_{i=1}^{s(P)-1}$, $\psi_{\inf}$ and $\psi_{\sup}$.
\end{itemize}
\end{algorithmic}
\end{minipage}}
\end{figure}

\paragraph{Notations}
In the following, we search iteratively for a time index corresponding to the minimal (resp. maximal) value of a subset of time series $A_k$. When there are several time indexes verifying this condition, we choose the smallest one. Hence, for $I$ a subset of $\{1,2, ..., N\}$, we define the following notations:
\[
\arg\min_{k\in I} A_k = \min \{j \in I | A_j = \min_{k \in I} A_k\} \\
\]
\[
\arg\max_{k\in I} A_k = \min \{j \in I | A_j = \max_{k \in I} A_k\}
\]

\paragraph{Step 1: Search for the first pulse.} Under the assumption that the maximum value of $A_k$ retrieved from a sufficiently large time window coincides with a pulse occurrence, the first pulse is detected by searching for the maximum value of $A_k$ for $k\in\{1,2,\dots,(2k_p)\}$ (where $k_p=\lfloor T_p/T_s \rfloor$). The choice of the time window size (twice nominal IPI) containing $2k_p$ samples is a trade-off between the risk of missing the first pulses (encountered when the window is too large) and the risk of false detection (encountered when the window is too small).
\begin{figure}[H]
\fbox{\begin{minipage}{\textwidth}
\begin{algorithmic}
\STATE Step 1: Find the first pulse occurrence index $p_1$ by searching for the maximum value of the sampled time series $A_k$ within the first $2 k_p$ data samples: 
\vspace{-0.2cm}
\[
  p_1:= \arg\max_{k\in\{1,2,\dots,(2k_p)\}} A_k
\]
\end{algorithmic}
\end{minipage}}
\end{figure}

\begin{figure}[htbp]
\centering
\includegraphics[width=17.5cm]{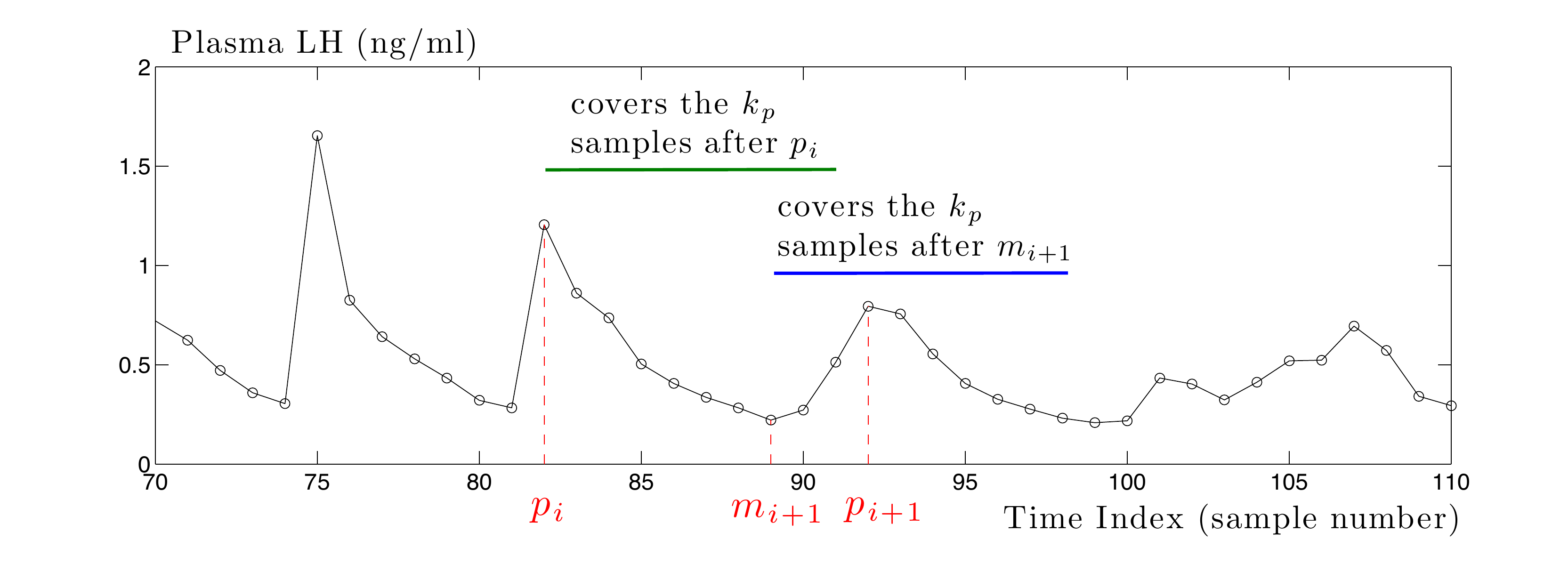}
\caption{
{\bf  One iteration of the forward research of pulses.} For a given value of $i$ in the iterative process (initialized with $i=1$), the algorithm searches for the index $m_{i+1}$ of the minimal sample from the $p_i$-th sample in the window defined by the nominal period $T_p$, i.e. among the $k_p$ samples (under the green segment) directly following the $p_i$-th sample. Then, the algorithm searches for the index $p_{i+1}$ of the maximal sample among the $k_p$ samples (under the blue segment) directly following the $m_{i+1}$-th sample. Index $p_{i+1}$ is stored in vector $P$ and the process is iterated with $i$ incremented by $1$ until the end of the time series.}
\label{jmin_jmax}
\end{figure}

\paragraph{Step 2: Search for the pulses following the first one.} After the detection of the first pulse located at $k=p_1$, the second one could be detected by searching for the maximum value of $A_k$ for $k\in\{(p_1+1),(p_1+2),\dots,(p_1+k_p+s)\}$, where some margin value $s$ should be chosen for the trade-off between over-detection and under-detection. To reduce the risks of misdetection, instead of searching for the maximum value of $A_k$ at this stage, the minimum value of $A_k$ for  $k\in\{(p_1+1),\dots,(p_1+k_p) \}$ is first looked for. Let $k=m_2$ be the location of this minimum, then the second pulse is searched for within the window $\{(m_2+1),\dots,(m_2+k_p) \} $. These searches are made within windows of size $k_p$ (see Figure \ref{jmin_jmax}). As the minimum and the maximum being searched for are expected to be inside the windows of size $k_p$, there is no need to choose a margin value. The following pulses are then searched for similarly.

\begin{figure}[H]
\fbox{\begin{minipage}{\textwidth}
\begin{algorithmic}
\STATE Step 2: Find the following pulse occurrences $p_i$ by alternatively searching for the minimum and maximum values of the sampled time series in a moving window covering  $k_p$ data samples:
\STATE \[ i:=1 \]
\vspace{-0.5cm}
\WHILE{$p_i+k_p \leq N$}
\STATE
    \[
      m_{i+1}:= \arg\min_{k\in\{p_i+1,p_i+2,\dots,p_i+k_p \}} A_k
    \]
    
    \IF{$m_{i+1} + k_p \leq N$}
    \STATE
      \[
         p_{i+1}:= \arg\max_{k\in\{m_{i+1}+1,m_{i+1}+2,\dots,m_{i+1}+k_p \}} A_k   
      \]
    \ENDIF
\STATE \[ i:=i+1 \]
\ENDWHILE
\vspace{-0.5cm}
\STATE  \begin{align*} 
              s_p &:= i-1 \\
                P &:= \{ p_1,\dots,p_{s_p}  \}
        \end{align*}
\end{algorithmic}
\end{minipage}}
\end{figure}

\newpage

\paragraph{Step 3: Remove too small peaks by various thresholding methods.} \quad

\paragraph{Step 3.1: Pulse height median-based thresholding.}
Several methods are used to remove small pulses possibly resulting from false detections. The first method is based on a threshold applied to the heights of the detected pulses. We recall that the height of the $i$-th pulse is defined as the difference between its amplitude $A_{p_i}$ and the lowest value of $A_k$ within the sampled time series, denoted by $\underline{A}$. If the height of a pulse is lower than the median of the heights of all the detected pulses multiplied by a ratio $\lambda_r$, then it is removed from the set of detected pulses.
\begin{figure}[H]
\fbox{\begin{minipage}{\textwidth}
\begin{algorithmic}        

\STATE Step 3.1. Pulse height median-based thresholding:
\begin{align*}
     & \underline{A}:= \min_{k\in{1,2,...,N}} A_k \\
     &   P:=   \{ p_i\in P:  A_{p_i}-\underline{A}  > \lambda_r (\text{median}(A_{p_1},\dots,A_{p_{s_p}}) - \underline{A} )\} \\
     & s_p:= s(P)
\end{align*}
\end{algorithmic}
\end{minipage}}
\end{figure}

\newpage 
 
\paragraph{Step 3.2: Local relative magnitude-based thresholding.}
The second method for removing false pulses is to compare the magnitude of each pulse with those of its neighbors. Let $B_1$ and $B_2$ be respectively  the minimum values of $A_k$ for $p_{i-1}<k<p_i$  and $p_i<k<p_{i+1}$. The local magnitude of the $i$-th pulse is defined as the geometric mean of $(A_{p_i} - B_1)$ and $(A_{p_i} - B_2)$. Let $B_0$ be the lowest value out of $B_1$ and $B_2$, The local magnitude of the $i$-th pulse is compared with the geometric mean of $(A_{p_{i-1}} - B_0)$ and $(A_{p_{i+1}} - B_0)$, relatively to a threshold $\lambda_r$.  The comparison of geometric means is equivalent to the comparison of arithmetic means in the logarithmic scale, so that it tends to favor smaller quantities.
\begin{figure}[H]
\fbox{\begin{minipage}{\textwidth}
\begin{algorithmic}  
\STATE Step 3.2. Local relative magnitude-based thresholding:
\FOR{$i=2,\dots,(s_p-1)$}
\vspace{-0.5cm}
\STATE 
\begin{align*}
     B_1 &:= \min_{k\in \{ (p_{i-1}+1),(p_{i-1}+2),\dots, (p_{i}-1)\}}  A_k \\
     B_2 &:= \min_{k\in \{ (p_{i}+1),(p_{i}+2),\dots, (p_{i+1}-1)\}}  A_k  \\
     B_0 &:= \min(B_1,B_2) \\
\end{align*}
\vspace{-0.5cm}
\IF{$(A_{p_i} - B_1)(A_{p_i} - B_2) < \lambda_r^2 (A_{p_{i-1}} - B_0)(A_{p_{i+1}} - B_0)$}
\STATE Erase $p_i$ from $P$ 
\ENDIF
\ENDFOR
\vspace{-0.5cm}
\STATE  \[ s_p:= s(P) \]
\end{algorithmic}
\end{minipage}}
\end{figure}

\paragraph{Step 3.3: Local magnitude absolute thresholding.} The third method for removing false pulses is based on some prior knowledge about the magnitudes of the pulses.  Let $B_1$ and $B_2$ be respectively  the minimum values of $A_k$ for $p_{i-1}<k<p_i$  and $p_i<k<p_{i+1}$. If the local magnitude of the $i$-th pulse (the geometric mean of $A_{p_i} - B_1$ and $A_{p_i} - B_2$) is smaller than a chosen threshold $\lambda_a$, then it is removed from the set of detected pulses.
\begin{figure}[H]
\fbox{\begin{minipage}{\textwidth}
\begin{algorithmic}  
\STATE Step 3.3. Local magnitude absolute thresholding:
\FOR{$i=2,\dots,(s_p-1)$}
\STATE 
\begin{align*}
     B_1 &:= \min_{k\in \{ (p_{i-1}+1),(p_{i-1}+2),\dots, (p_{i}-1)\}}  A_k \\
     B_2 &:= \min_{k\in \{ (p_{i}+1),(p_{i}+2),\dots, (p_{i+1}-1)\}}  A_k
\end{align*}
\IF{$(A_{p_i} - B_1)(A_{p_i} - B_2) < \lambda_a^2$}
\STATE Erase $p_i$ from $P$ 
\ENDIF
\ENDFOR
\STATE \[  s_p:= s(P) \]
\end{algorithmic}
\end{minipage}}
\end{figure}

\newpage

\paragraph{Step 4: Retrieve missed pulses.}
To retrieve possible missed pulses, the data points lying between each pair of detected pulses are examined. For each sample index $j$ lying between the pulse occurrences $p_i$ and $p_{i+1}$, the relative height of the sample $A_j$ is defined as the geometric mean of $A_j-B_1$ and $A_j-B_2$, where $B_1$ and $B_2$ are respectively the lowest values of $A_k$ for $p_i<k\leq j$ and for $j\leq k <p_{i+1}$. The sample exhibiting the maximum relative height between $p_i$ and $p_{i+1}$ is compared with the geometric mean of $A(p_i)-B_0$ and $A(p_{i+1})-B_0 $, where $B_0$ is the lowest sample value between $p_i$ and $p_{i+1}$. If the comparison with respect to the threshold  $\lambda_r$ is conclusive, a new pulse is added to the set of detected pulses. This process is repeated three times to deal with the case where several pulses might have been missed between two consecutive previously detected pulses..
\begin{figure}[H]
\fbox{\begin{minipage}{\textwidth}
\begin{algorithmic}  
\STATE Step 4. Retrieve missed pulses.
\FOR{$k=1,2,3$}
\STATE $J=\emptyset$ \\
\FOR{$i=1,\dots,(s_p-1)$}
\IF{$p_i+3<p_{i+1}$}
   \FOR{$j=(p_i+2),\dots,(p_{i+1}-2)$}
     \STATE 
       \begin{align*} B_1 &:= \min_{k\in\{ (p_i+1),\dots,j \}} A_k \\
                      B_2 &:= \min_{k\in\{ j,\dots,(p_{i+1}-1) \}}  A_k\\
                      C_j &:= (A_j-B_1)(A_j-B_2) 
       \end{align*}
   \ENDFOR
   \STATE
   \begin{align*}
    & j_{\text{max}}:= \arg\max_{j\in\{ (p_i+2),\dots,(p_{i+1}-2)\}}  C_j \\
    & B_0:= \min(B_1, B_2)
   \end{align*}
   \IF{$C_{j_{\text{max}}} > \lambda_r^2 (A(p_i)-B_0)(A(p_{i+1})-B_0)$}
     \STATE Insert $j_{\text{max}}$ into $J$
   \ENDIF
\ENDIF
\ENDFOR
\STATE Insert $J$ into $P$ \\
\STATE $s_p=s(P)$
\ENDFOR 
\end{algorithmic}
\end{minipage}}
\end{figure}

\newpage
 
\paragraph{Step 5: Removal of 3-point peaks.}
When the rhythm of the hormonal pulses accelerates, the IPI may approach the sampling period $T_s$, to the point that it becomes impossible to reliably detect pulses. In such situations, the sampled time series may exhibit local maxima supported by 3 consecutive data samples. To avoid false detections, pulses detected upon 3 points are identified and possibly removed. For a pulse detected at $k=p_i$, if $p_i-1$ and $p_i+1$ are local minima ``sharp'' enough, in the sense that $A(p_i-2)-A(p_i-1)$ and $A(p_i+2)-A(p_i+1)$ are large enough compared to $A(p_i)-A(p_i-1)$ and $A(p_i)-A(p_i+1)$, then $A_{p_i}$ is considered as a peak detected upon 3 points (referred to as ``3-point peak'') and is removed from the set of detected pulses.
\begin{figure}[H]
\fbox{\begin{minipage}{\textwidth}
\begin{algorithmic}
\STATE Step 5. Remove pulses detected upon 3 data points:
\vspace{-0.4cm} 
\STATE \[  s_p:= s(P) \]
\vspace{-0.5cm} 
\FOR{$i=1,\dots,s_p$}
\vspace{0.2cm}
   \IF{$
      \left\{
      \begin{array}{llll}
        &  p_i \geq 3, &
      \text{and } & p_i \leq N-2, \\
      \text{and } & A(p_i-2)>A(p_i-1), &
      \text{and } & A(p_i)>A(p_i-1), \\
      \text{and } & A(p_i)>A(p_i+1), &
      \text{and } & A(p_i+2)>A(p_i+1)
      \end{array}
      \right\}
      $}
     \vspace{0.4cm}
     \STATE \[R:=\frac{ \{[A(p_i-2)-A(p_i-1)] + [A(p_i+2)-A(p_i+1)] \}/2}{\sqrt{[A(p_i)-A(p_i-1)][A(p_i)-A(p_i+1)]}}\]
     \IF{$R\geq \lambda_{3p}$}
        \STATE Erase $p_i$  from $P$
     \ENDIF
   \ENDIF
\ENDFOR
\end{algorithmic}
\end{minipage}}
\end{figure}

\paragraph{Step 6: IPI sequence and tunnel construction.}
Under some regularity assumptions of the IPIs, the detected pulses leading to IPI outliers should be corrected. For this purpose, a cubic polynomial $\theta_0 + \theta_1 x+ \theta_2 x^2 + \theta_3 x^3$ is first fitted to the IPIs where $x$ is equal to the shifted pulse index $i$. Let $(\hat\theta_0, ... , \hat\theta_3)$ be the value of $(\theta_0, ..., \theta_3)$ fitted to the segmented IPIs and $\phi(x)=\hat\theta_0 + \hat\theta_1 x+ \hat\theta_2 x^2 + \hat\theta_3 x^3$. The lower and upper edges of the IPI tunnel are defined respectively by
$\phi \left(x+\frac{s_p-1}{2}\right)(1-\alpha)$ and $\phi \left(x+\frac{s_p-1}{2}\right)(1-\beta)$. The centering of the $i$-indexes (leading to $x_i$) is for the purpose of a better numerical accuracy.
\begin{figure}[H]
\fbox{\begin{minipage}{\textwidth}
\begin{algorithmic}
\STATE Step 6. IPI sequence and tunnel construction: 
\vspace{-0.4cm}
\STATE \[  s_p:= s(P) \] 
\vspace{-0.5cm}
\FOR{$i=1,\dots,s_p-1$}
\vspace{-0.7cm}
   \STATE \begin{align*} 
            q_i &:= p_{i+1}-p_i \\
            x_i &:= i - \frac{s_p-1}{2}
          \end{align*}
\vspace{-0.7cm} 
\ENDFOR
\vspace{-0.9cm}
\STATE \begin{align*}
        (\hat\theta_0, \hat\theta_1, \hat\theta_2, \hat\theta_3) =\arg \min_{(\theta_0, ... , \theta_3) \in\mathbb{R}^4}\sum_{i=1}^{s_p-1} \left[ (\theta_0 + \theta_1 x_i+ \theta_2 x_i^2 + \theta_3 x_i^3) - q_i  \right]^2     
       \end{align*}
\end{algorithmic}
\end{minipage}}
\end{figure}

\subsection*{Experimental data provided to the algorithm}

We have run the algorithm on either experimental or synthetic LH time series. Experimental time series included nineteen ewes, distributed over two different protocols. All procedures were approved by the ``Direction D\'epartementale des Services V\'et\'erinaires d'Indre-et-Loire'' (approval number C37-175-2)  for the agricultural and scientific research agencies INRA (French National Institute for Agricultural Research) and CNRS (French National Center for Scientific Research),  and conducted in accordance with the Guide for the Care and Use of Agricultural Animals in Research and Teaching. Blood samples from a first group of ten estrus-synchronized ewes (Lacaune breed, \cite{Drouilhet_10}) were collected via jugular venous cannula every $10$ min for a period of $24$ h during the follicular phase. A second group of nine ovariectomized Ile-de-France ewes were collected during anestrus season for blood sampling every 10 min over a period of 15h. Ewes received an agonist of somatostatin type 2 receptor via intracerebroventricular injection between 5 and 10h after sampling start (Courtesy of A. Caraty, unpublished data). All blood samples were collected into heparinized tubes and then centrifuged for $20$ min at $400$ g. Plasma was stored at $-20$\degre C until hormone assays \cite{Drouilhet_10}.

\section*{Results}

The output of our algorithm consists of the IPI series, providing the number of detected peaks with respect to the time series indexes. Moreover, the IPI tunnel has been used on the LH time series according to the assumption of regularity in their frequency modulation. The sampling period and the absolute magnitude threshold $\lambda _a$, corresponding to the minimal detectable concentration, are provided by the protocol specifications. The default set, proposed in Table \ref{tab:ParamVal}, has been used in all the cases: nominal period, $T_p$, equal to $40$ min, relative magnitude threshold, $\lambda _r$, equal to $0.2$, ($20$\% of the geometric mean of the neighbors), both the lower and upper bound of the IPI tunnel equal to $0.6$. The choice of the default parameter set is explained.

\subsection*{The InterPulse Interval (IPI) series}

\begin{figure}[!ht]
\centering
\includegraphics[width=15cm]{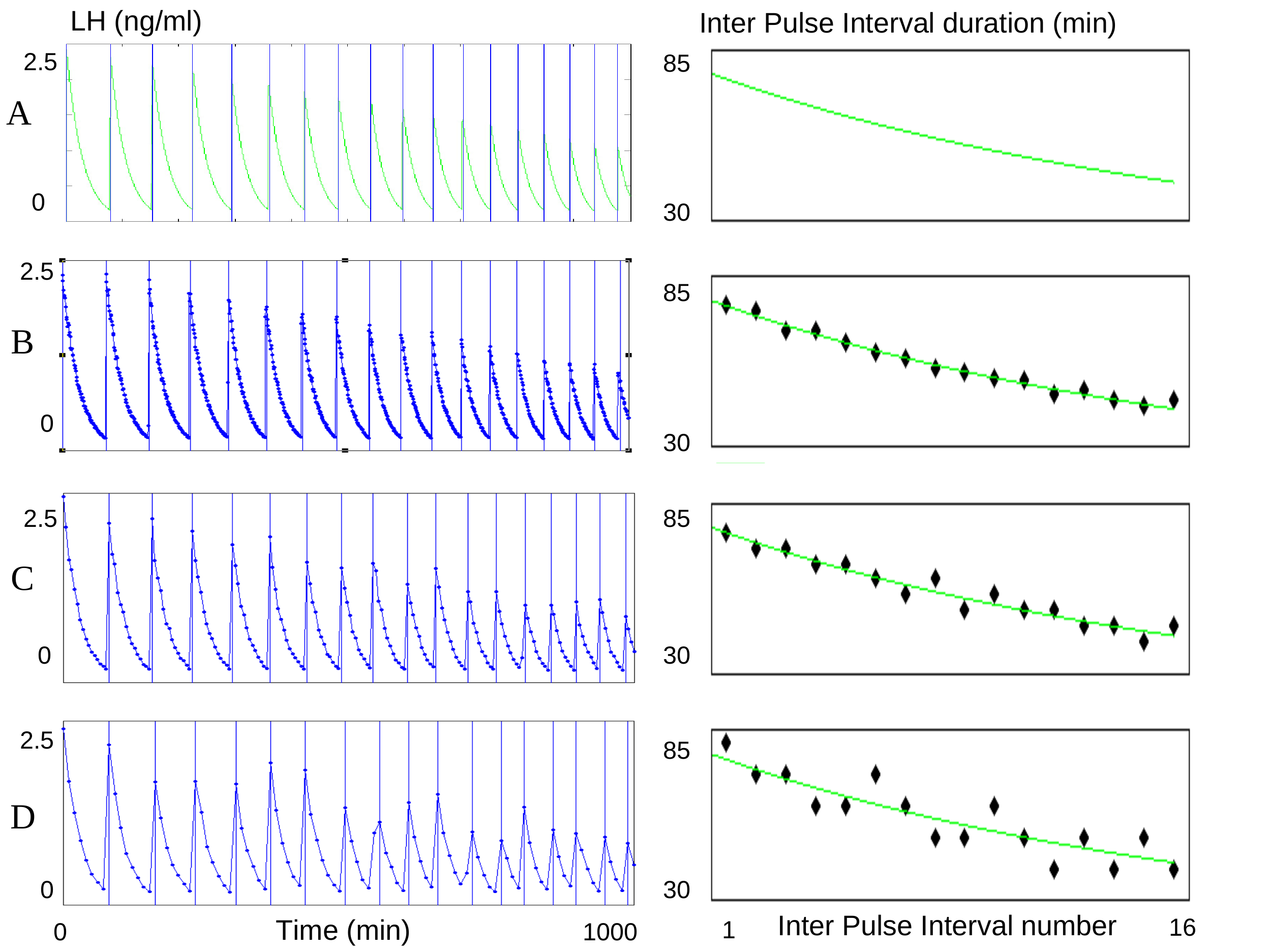}
\caption{
 {\bf  IPI series from synthetic LH time series with different sampling frequencies.}
Left panels: LH plasma level time series retrieved over $1000$ min. Vertical lines correspond to the pulse occurrences. Panel A: theoretical plasma level, corresponding to a continuous monitoring. The pulse frequency increases, whereas the pulse amplitude decreases along time.  Panels B, C and D: sampled series, with a respective sampling period of $1$, $5$ or $10$ min.  Stars on the time series correspond to sampled points. The first sampling time occurs at the first minute of the simulation ($r= 1$ min), variability in the sampling times is set to $15$\% of the sampling period ($f=0.15$ min for B, $f=0.75$ min for C and $f=1.5$ min for D) and the assay variability is set to $b=5$\%. \newline
Right panels:  resulting IPI series, indexed by the number of the pulse occurrence (each IPI is represented by a black diamond). The theoretical IPI series is the continuous green curve, superimposed on the IPI series obtained after sampling. In any case, there are 16 detected peaks.
}
\label{IPI_simul}
\end{figure}

\begin{figure}[!ht]
\centering
\includegraphics[width=17cm]{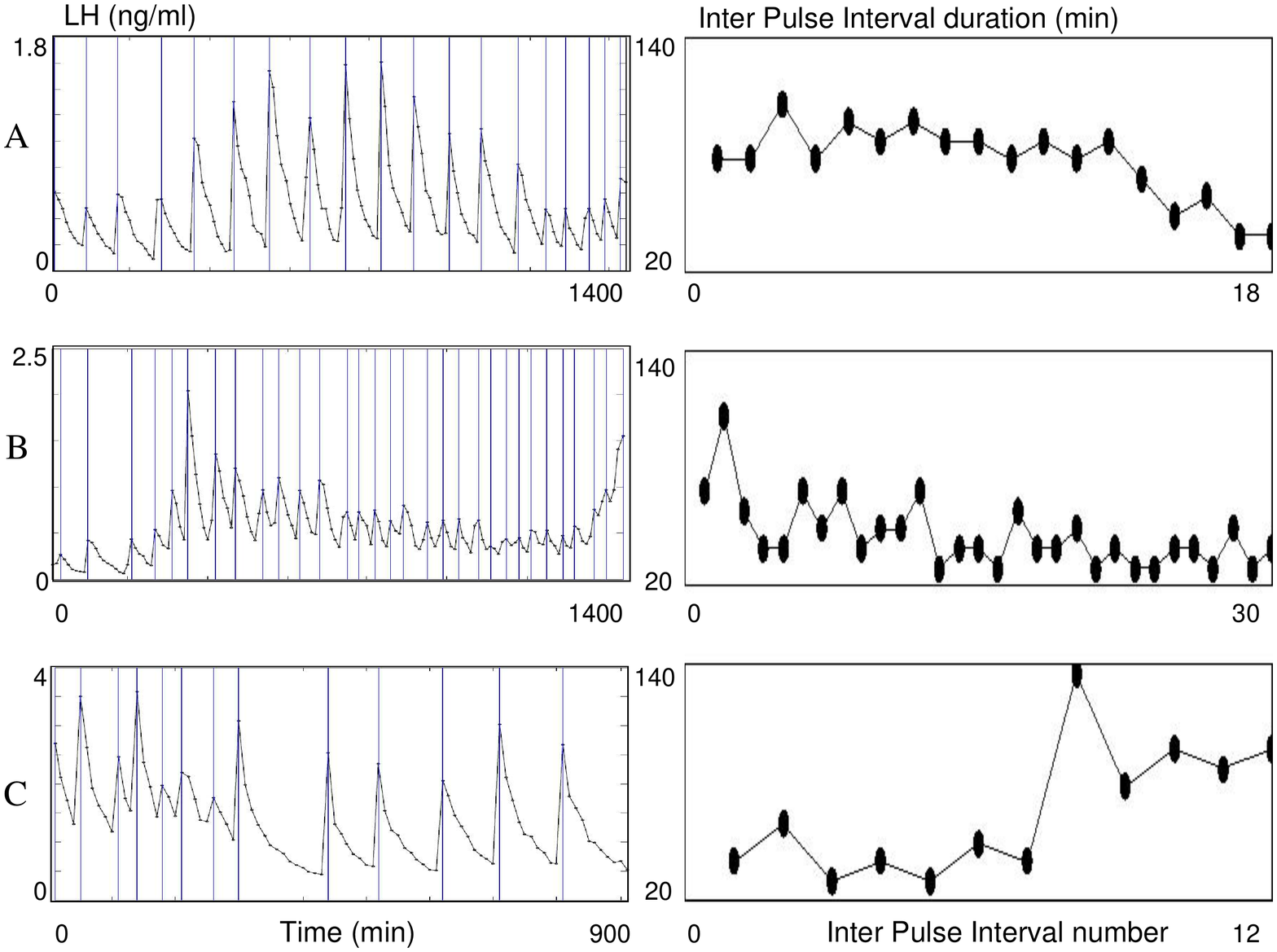}
\caption{
 {\bf  IPI series from experimental LH series with different pulsatile rhythms.}
Left panels:  LH plasma level time series with a $10$ min. sampling period. Vertical lines correspond to the pulse occurrences. Stars on the time series correspond to sampled points. Right panels: resulting IPI series (black diamond), indexed by the number of the pulse occurrence. Panel A: stable rhythm with final acceleration. Panel B: progressive acceleration. Panel C: fast deceleration in the second half of the series.
}
\label{IPI_reels}
\end{figure}

Figures \ref{IPI_simul} and \ref{IPI_reels} correspond respectively to LH synthetic and experimental series. The left panels display the LH plasma level time series.Vertical lines correspond to the pulse occurrences. Stars on the time series correspond to the points of sampled measures. The right panels display the resulting IPI series, indexed by the number of the pulse occurrence (each IPI is represented by a black diamond).


On the synthetic LH series (Figure \ref{IPI_simul}), the left panels display the following cases. Panel A represents a theoretical plasma level series, that would be retrieved in case of continuous monitoring. The pulse frequency increases whereas the pulses amplitude decreases along time. Panels B, C and D are the corresponding sampled series with a respective sampling period of $1$, $5$ or $10$ min. The first sampling time occurs at the first minute of the simulation ($r= 1$ min), and there is variability both in the sampling times ($f$ is equal to $15$\% of the sampling period, corresponding to $0.15$ min for B, $0.75$ min for C and $1.5$ min for D) and the assays ($b= 5$\%). On the right panels, the theoretical IPI series are represented by a continuous green curve, superimposed on the IPI series of measured LH series (black diamonds). Comparisons between the successive panels allow us to assess the influence of the sampling period on the IPI series. The number of detected peaks ($16$) is the same, and the patterns of regular acceleration are identical, whatever the sampling period is.\\
Moreover, the discretization of the initially continuous signal induces delays in the time occurrence of pulses; the maximal delay corresponds to the sampling period. The higher the sampling period is, the closer the measured IPI series are to the theoretical ones.
On the experimental LH series (Figure \ref{IPI_reels}), three different pulsatile rhythms are displayed, with a $10$ min sampling period. Panel A illustrates the case of a stable rhythm with a final acceleration resulting in IPI shortening in the last third of the series.  Panel B illustrates the case of a progressive acceleration resulting in a progressive shortening of the IPIs. Panel C illustrates the case of a fast deceleration in the second half of the series resulting in increased IPIs.  As the IPI series give information both on the number of detected peaks and the rhythm evolution, they are particularly useful for comparing the rhythmicity of different series of the same duration (for instance, B is almost twice as fast as A).

\subsection*{The IPI tunnel}
 
Due to the assumption of regularity in the frequency modulation of the LH time series, the IPI tunnel has been used to point out situations where there may be a lack of detection or an over-detection of pulses in the time series. In such situations, we can try to explain the detection error and propose possible corrections. On the opposite, the IPI tunnel can detect genuine long or short IPIs, and be used as a tool for analyzing sudden frequency breaks or accelerations  in pulsatile rhythms. In all examples the sampling period was equal to $10$ min.

\begin{figure}[!ht]
\centering
\includegraphics[width=16.3cm]{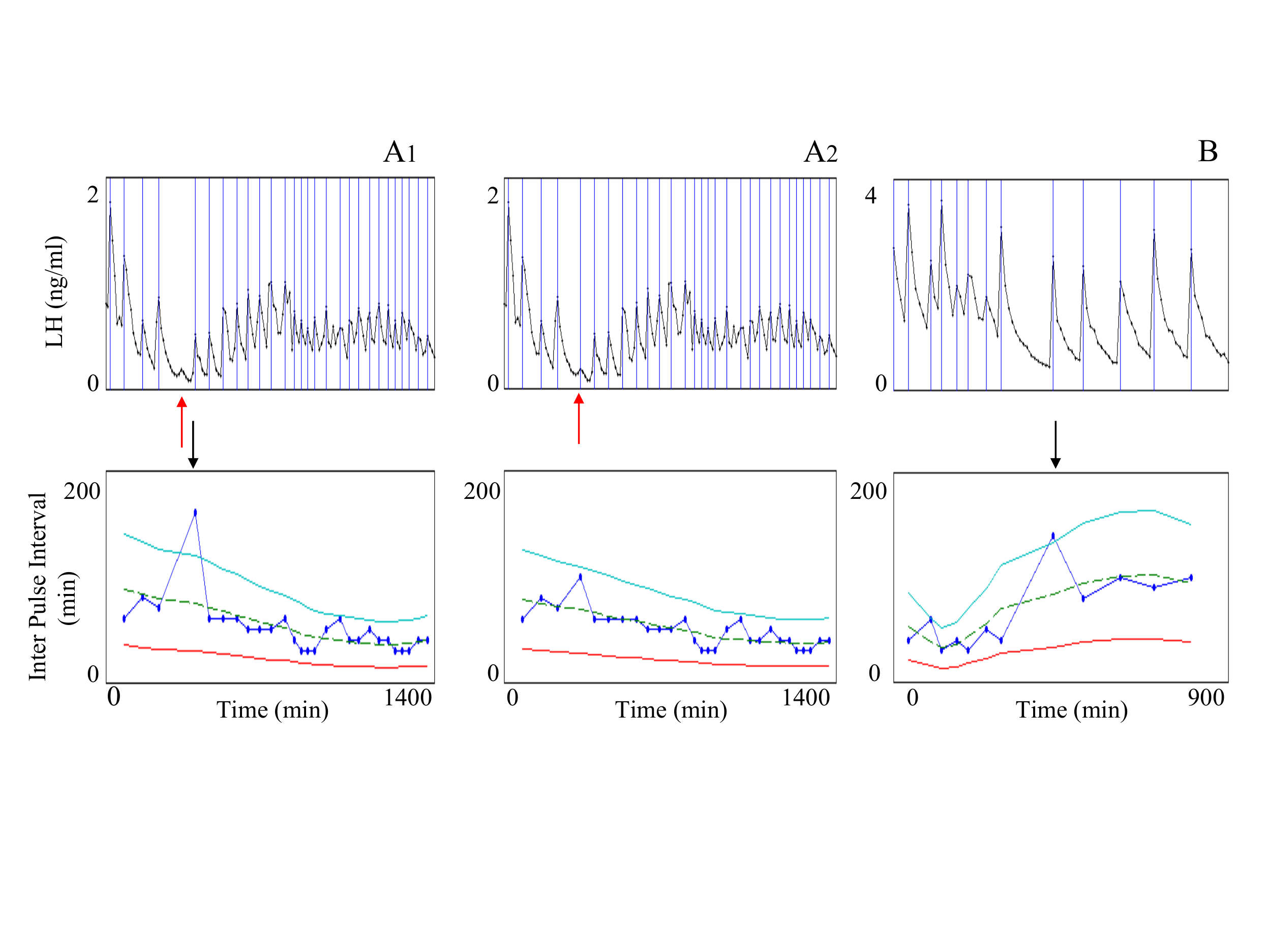}
\vspace{-2.5cm.pdf}
\caption{{\bf  IPI outliers lying above the upper bound of the tunnel. Example of correction by  decreasing the value of the relative magnitude threshold, $\lambda _r$.}
Top panels: two  experimental LH plasma time series, A (panels A1 and A2) and B (panel B).  Vertical lines correspond to pulse occurrences. Bottom panels: resulting IPI series indexed by time. Each IPI value (blue point) corresponds to the time elapsed between the current detected pulse and the previous one. Dashed line: moving cubic function fitting the values of the IPI series. Solid lines: lower ($\alpha$-dependent function $\psi_{\inf}$) and upper ($\beta$-dependent function $\psi_{\sup}$) edges of the tunnel ($\alpha = \beta = 0.6$). Black arrows: occurrences of the outliers. Case A: outlier due to a lack of detection (missed pulse designed by a red arrow); panel A1: initial IPI series with $\lambda _r=0.2$ (default value); A2: corrected IPI series; with $\lambda _r=0.1$. Case B: genuine long IPI.}
\label{outliers_sup}
\end{figure}

\begin{figure}[H]
\centering
\includegraphics[width=16.3cm]{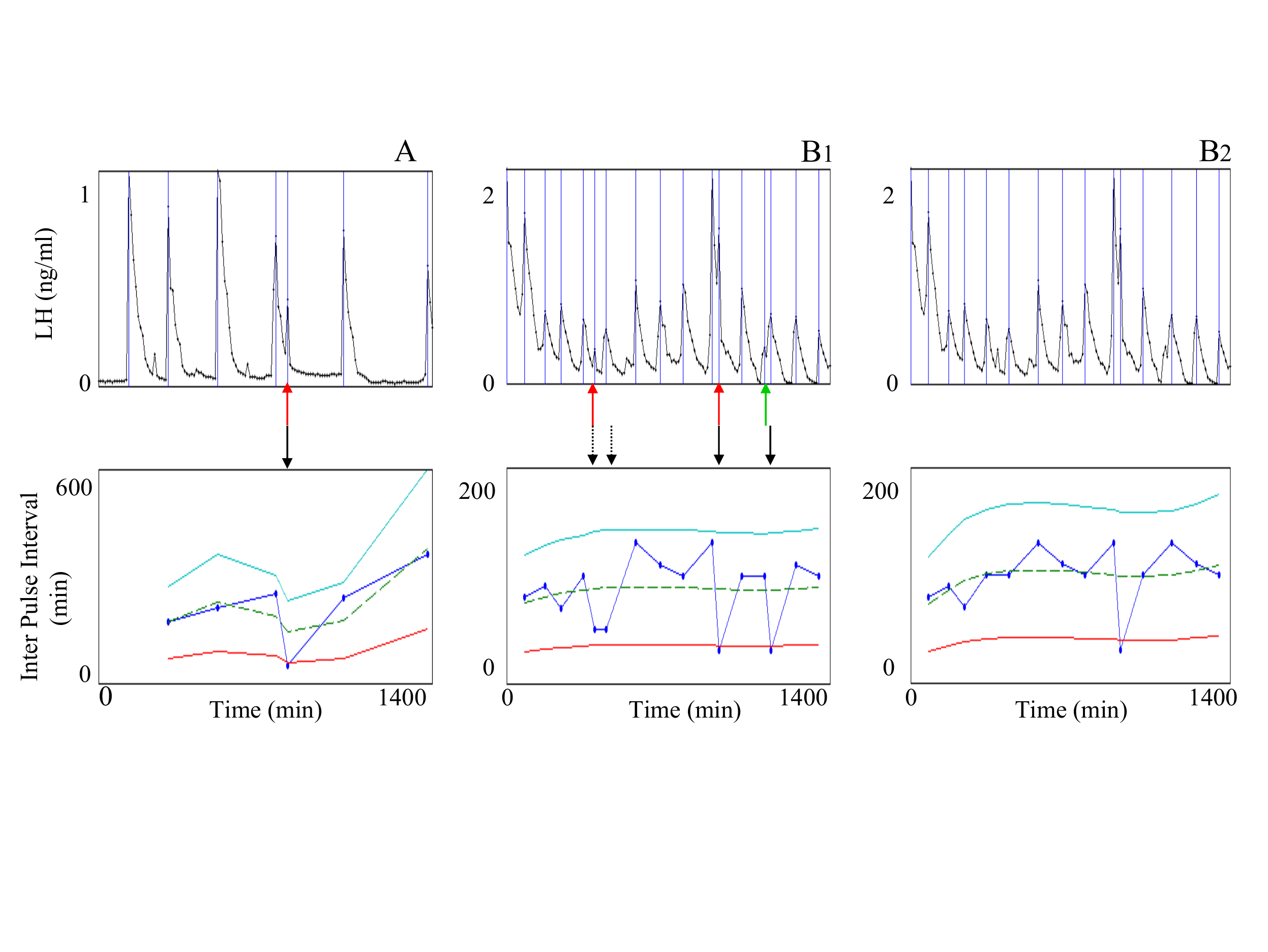}
\vspace{-2.5cm}
\caption{
 {\bf  IPI outliers lying below the lower bound of the tunnel. Example of correction by increasing the value of the relative magnitude threshold, $\lambda _r$.}
Top panels: two experimental LH plasma time series, A (panel A) and B (panels B1 and B2). Bottom panels: resulting IPI series indexed by time. The vertical lines, the blue points, the dashed line and the solid lines represent the same objects as in Figure \ref{outliers_sup}. 
Panels A and B1: initial IPI time series. Solid black arrows: occurrence of clear outliers. Dashed black arrows: occurrence of  IPIs that can be associated with over-detected peaks in the LH series although they remain above the lower bound of the tunnel. Red arrows: peaks lying on the middle of the descending phase of the preceding pulse. Green arrow: peak lying on the middle of the ascending phase of the following pulse. 
Panel B2: B corrected IPI series after increasing the relative magnitude threshold $\lambda_r$ to 0.45. Two of the three false peaks have been discarded.
}
\label{outliers_inf}
\end{figure}

Figures \ref{outliers_sup}  and \ref{outliers_inf}  display respectively apparent lacks of detection or over-detections in experimental time series. The top panels display the LH plasma time series and vertical lines correspond to pulse occurrences. The bottom panels display the resulting IPI series, indexed by the pulse time rather than the pulse number in order to keep the same reference time in both the LH and IPI series. Each IPI value (marked by a blue point) corresponds to the time elapsed between the current detected pulse and the previous one. The IPI series are displayed together with the three curves delimiting the tunnel: the dashed line represents the moving cubic function fitting the values of the IPI series, the solid lines represent the lower ($\psi_{\alpha}$) and lower ($\psi_{\beta}$) bounds of the tunnel.

Figure \ref{outliers_sup} displays two cases of apparent lack of detection, where the IPI outliers lie above the upper bound. 
In case A, the outlier appears at minute $400$ (black arrow, panel A1). Going back and forth between the IPI series and the LH time series allows us to favor the hypothesis of a lack of detection. Indeed, if we take into account the small amplitude pulse occurring at minute $340$ (red arrow, panel A1),  the exceedingly large IPI can be distributed over two consecutive IPIs of $100$ and $60$ min, whose duration are compatible with the local tunnel size (local upper bound of $123$ min, local lower bound of $30.5$ min), hence with the regularity assumption.

\vspace{-1.5cm}

A first correction step consists in decreasing the relative magnitude threshold ($\lambda _r$), set to the default  value of $0.2$, in such a way that the missed pulse can be recovered without adding false detections. Panel A2 illustrates the result of the correction: the missed pulse occurring at minute $340$ was recovered (red arrow) after decreasing the value of parameter $\lambda _r$ to $0.1$.
In case B, the outlier appears at minute $440$ (black arrow, panel B). Going back and forth between the IPI series and the LH time series allows us to favor the hypothesis of a genuine long IPI. There is not only no visible pulse after the preceding detected pulse but the rhythm also remains slow after the long IPI.

Figure \ref{outliers_inf} displays some cases of over-detection, where the IPI outliers lie below the lower bound, in two LH series A and B. The outlier occurrences are indicated by solid black arrows.  It is worth noticing that the two IPIs indicated by dashed arrows in case B (panel B1)  are not classified as outliers with the default set parameters ($\alpha = \beta = 0.6$) but they are close enough to the tunnel lower bound to draw the user's attention. In this example, we considered them as IPI outliers. In both cases, the patterns of some detected peaks suggest false detections possibly imputable to measurement conditions.

\vspace{-1.8cm}

The IPI outliers in cases A (panel A) and B (panel B1) are either due to additional peaks lying on the middle of the descending phase of the preceding pulse (red arrows) or  to a peak lying on the middle of the ascending phase of the following pulse (green arrow). A first correction step consists in increasing the relative magnitude threshold ($\lambda _r$) in such a way that the peak can be discarded without eliminating true detections. Increasing $\lambda _r$ to $0.45$, allows us to discard two false peaks in case B (panel B2).   Nevertheless, no change in $\lambda _r$ can get rid of the IPI outlier in case A nor of the third IPI outlier in case B without discarding genuine pulses at the same time. It appears that the amplitude of the peaks underlying the IPI outlier is too close to that of their neighbors. It will be up to the user to evaluate the influence of such IPIs on the characteristics of the series and to follow the more appropriate strategy.

\begin{figure}[!ht]
\centering
\includegraphics[width=15.7cm]{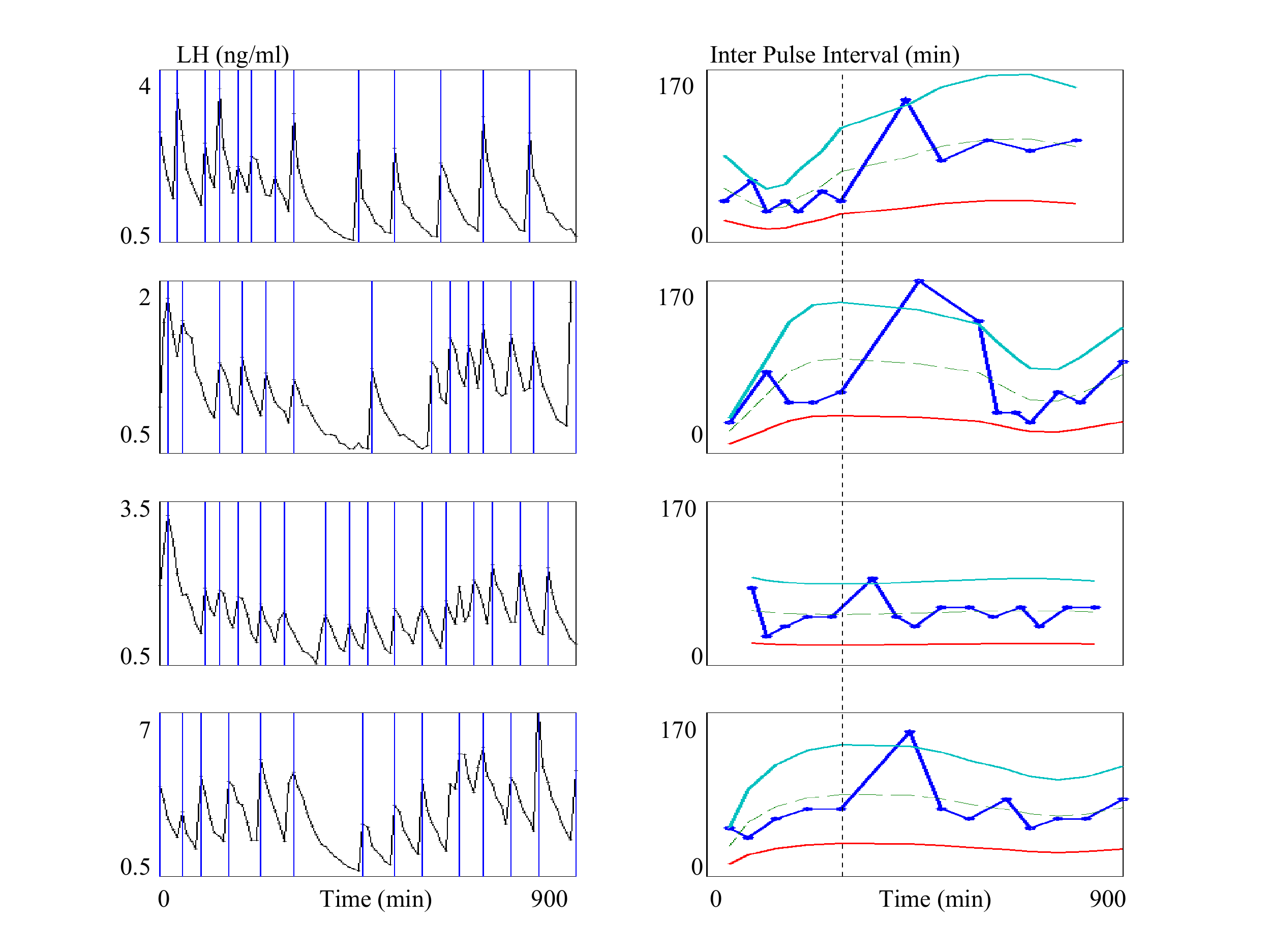}
\caption{
 {\bf  IPI-based study of the synchronization between LH series.}
The four LH series are retrieved from ewes subject to an experimental protocol inducing a steep decrease in the pulse frequency.
Left panels: experimental LH plasma time series; vertical lines correspond to pulse occurrences. Right panels: resulting IPI series indexed by time. The vertical lines, the dashed line and the solid lines represent the same objects as in Figure \ref{outliers_sup}. Vertical dashed line: break in the dynamics of the IPI series corresponding to the last IPI preceding the outlier.
}
\label{ruptures}
\end{figure}

Figure \ref{ruptures} shows how the IPI tunnel can be used for detecting sudden frequency breaks. It displays four different LH time series (left panels) retrieved from ewes subject to an experimental protocol inducing a steep decrease in the pulse frequency.

\vspace{-0.3cm}

IPI series (right panels) are indexed by the pulse time in order to keep the same reference time when studying variability between series. The break in the dynamics of the IPI series corresponds to the occurrence of the last IPI preceding the outlier lying above the upper bound. Moreover, identifying its precise location enables us to study the synchronization between the different time series, as pointed out by the vertical dashed line.

\subsection*{User parameters: choice of a default set and robustness evaluation}

Among the algorithm parameters (Table \ref{tab:ParamVal}), two are provided directly by the protocol specifications: the sampling period $T_s$ and the absolute magnitude threshold $\lambda _a$, corresponding to the assay detection threshold (minimal LH level that can be reliably measured). For $T_s$, an upper bound equal to $10$ min is recommended (see the explanation below). Default values have been fixed for the other five parameters of Table \ref{tab:ParamVal}: nominal period ($T_p$), relative magnitude threshold ($\lambda _r$), 3-point peak threshold ($\lambda_{3p}$) and the lower ($\alpha $) and upper ($\beta $) bounds of the IPI tunnel. The choice was based on the observation of LH time series retrieved in 19 ewes distributed over two different protocols.

\begin{table}[!ht]
\centering
\begin{tabular}{|c|c|c|c|c|c|}
\hline
 & Name & Notation & Default set & A & B \\
\hline
Time series characteristics & Sampling frequency                 & $T_s$              &  & X & X \\
                                                 & Absolute magnitude treshold & $\lambda_a$  &  & X & X \\
\hline
Parameters  & Nominal period                            & $T_p$ & 40  & X & X \\
                       & Relative magnitude threshold  & $\lambda_r$ & 0.2  & X & X \\
                       & 3-point peak threshold              & $\lambda_{3p}$ & 0.1  & X &  \\
                       & Frequency regularity of pulsatility & $\alpha, \beta$ & 0.6  & X &  \\
\hline
\end{tabular}
\caption{
{\bf Algorithm parameters and default set adapted to LH time series.} Default set for $T_s$ and $\lambda _a$ are provided by the experimental protocol specifications.
A: needed parameters for time series characterized by  a pulsatile pattern, an asymmetric shape of pulses, some regularity in the time evolution of the pulse frequency (LH, GH, Insulin).
B: needed parameters for time series with no underlying assumption of symmetric shape of the peaks or frequency regularity (intracellular calcium).}
\label{tab:ParamVal}
 \end{table}

\paragraph{Nominal period ($T_p$)}
$T_p$ is chosen so as to favor the maximal number of correct detections, especially at the beginning of time series. According to the existence of high frequency series with short IPIs, there is a risk to detect two consecutive pulses in the same window if it is too large. $T_p$ has been set to $40$ min, which is a value close to the minimal observed period. The number of missed pulses is equal to $0$ for $T_p$ ranging from $40$ to $70$ min; it increases up to $2$ for $T_p=80$ min.  Even in that latest case, the number of missed pulses is small enough to guarantee the algorithm robustness with respect to this parameter.

\paragraph{Relative magnitude threshold ($\lambda _r$)}
Compared to the performances of the algorithm with the $\lambda _r$ default value ($0.2$), a decrease down to $0.1$ or an increase up to $0.3$ increase the total number of outliers (either false detections or missed detections). Moreover, this parameter can be adapted in order to correct outliers lying outside the tunnel, as previously seen.

\paragraph{3-point peak threshold ($\lambda_{3p}$)}
The 3-point peak threshold is introduced to prevent non-asym\-metric pulses from being detected. Figure \ref{3-point-peak} illustrates the identification of a 3-point peak pattern. $\lambda_{3p}$ is the ratio between the arithmetic mean of the amplitude of the neighboring points of rank 2 (green lines) and the geometric mean of the amplitude of the immediate neighbors (red lines). For instance, based on that criterion, the peak selected on panel A (dashed vertical line), is identified as a genuine 3-point peak. On the contrary, the peak selected on Figure \ref{3-point-peak}, panel B (solid vertical line) is not identified as a 3-point peak, since it belongs to a genuine, asymmetric LH pulse with an exponential decrease, albeit locally noised (the LH level corresponding to the latter neighbor of rank 2 is a little higher than expected for a smooth exponential decrease). Over 15 analyzed potential 3-point peaks, the minimal $\lambda_{3p}$ value corresponding to a genuine 3-point peak was equal to $0.15$, while the maximal $\lambda_{3p}$ value corresponding to a locally noised, genuine asymmetric LH pulse was equal to $0.07$, so that there was no overlapping. Consistently, we chose a default set value ($0.1$) lying within the $[0.07,0.15]$ range. This parameter is embedded within the algorithm, since there is no reason for the user to modify it. Indeed, the 3-point peaks rely on endocrinological considerations and take into account, either directly or indirectly, the typical duration of a LH pulse (around $30$ min) and its asymmetric shape, that can be reconstructed from time series sampled at least every $10$ min.

\paragraph{Lower ($\alpha$) and upper ($\beta$) bounds of the IPI tunnel} For the lower bound $\alpha $, the $0.6$ default value does not lead to any IPI outlier, whereas a value of $0.5$ leads to classify two genuine pulses as over-detected pulses, and a value of $0.4$ leads to classify twelve genuine pulses as over-detected pulses (among more than $300$ pulses). For the upper bound $\beta$, the $0.6$ default value allows one to identify every genuine outliers, without generating false outliers, but changing the value to $0.5$ led to false additional under-detected pulses. However, the use of the tunnel bounds is mainly to draw the user's attention to possible events of interest, and there is no direct sensitivity of the algorithm to their precise values.

\section*{Discussion}

For hormonal time series as those studied in this article, an important particularity is the fact that the signals are clearly subsampled due to the invasive nature of the sample collection procedure. In this case the classical filtering methods cannot be successful. Specific methods have to be developed to overcome this difficulty. There are two main approaches to study time series of pulsatile hormones. One consists in trying to detect, as accurately as possible, the pulse peaks, considered as discrete events \cite{urban_88}, while the other is based on deconvolution principles and intend to reconstruct the underlying secretion process \cite{veldhuis_90}. The deconvolution approach might seem more attractive, since it is susceptible to provide rich information on the hormonal signal, but it is hampered by the lack of validation, since information on the ``true'' signal is almost never available and cannot be directly compared to the reconstructed signal. Our own algorithm clearly belongs to the category of discrete peak detectors. Whereas, to our knowledge, the other available algorithms rely only on local and semi-local amplitude criteria, our algorithm combines local (on the data point level), semi-local (on the level of -possibly moving- windows of consecutive data points), and global (on the whole series level) amplitude criteria, with other criteria accounting for the pulse duration and the relative regularity in the pulse frequency modulation. Hence, this is a multi-scale and multi-criteria algorithm based on a dynamical selection process of the peaks.

To design our algorithm, the first idea was to locate significant local maxima in the processed time series, guided by the nominal IPI value so that the detected pulses have a reasonable rhythm. As the nominal IPI value may be too large or too small for each actual IPI in the processed signal, more steps were added to retrieve missed pulses and to remove false pulses. These extra steps are mainly based on the height of each pulse candidate and its magnitude with respect to the relative magnitude between neighbor pulses. The main advantage of these criteria is their weak dependence on the baseline that is most of the time non-stationary in typical hormonal time series.

An original feature of our work is to combine mathematical modeling with signal processing.  We have used synthetic time series, generated by a simple dynamical model, to illustrate the fundamental concepts underlying our algorithmic choices, as well as to assess the robustness of the model outputs with respect to the sampling rate and different sources of variability. Introducing uncertainty on both the measured LH level (to mimic assay variability) and time of measurement (to mimic possible hidden variability in the sampling chronology) allowed us to check that the ability to detect the right number of events was not affected by the noise.

For a given time series, the outputs of the algorithm consists of a corresponding series of detected IPI, structurally expressed as multiples of the sampling period. Hence, the algorithm provides information on the evolution of the frequency regime along the series, which is essential for studying the control of frequency encoding in endocrine systems.

We have run the algorithm on two different sets of experimental time series collected in sheep. Since they have a large body size and a much longer ovarian cycle compared to rodents, domestic species such as the ovine species are more suited to longitudinal endocrine studies and  their reproductive physiology is much closer to human reproductive physiology. We gave several instances showing that the algorithm is able to adapt to different patterns of frequency modulation (more or less rapid acceleration or deceleration) and also to detect breaks in the IPI rhythm. We then explained how one can make use of the IPI tunnel to discriminate outlier pulses from genuine pulses corresponding to a locally marked change in the frequency regime. On the whole, these results have shown that the algorithm can be employed to study and understand the frequency encoding of hormonal signals. To put the algorithm at the disposal of the user not familiar with computer programs, we intend to develop a user-friendly interface to make our software easily available and ready for use. The aim of this tool is to provide as much aid to decision as possible to the users together with guaranteeing full understanding on the detection process and the effect of the parameter values on the output.

In addition to the time series, the sampling period $T_s$ and the assay detection threshold $\lambda_a$ provided directly by the protocol specifications, there are only $5$ parameters to be set by the user: $T_p$, $\lambda_r $, $\lambda_{3p}$ for the pulse detection itself, $\alpha $ and $\beta $ for the definition of the IPI tunnel edges. A default set of parameter values is proposed in the case of LH time series (Table \ref{tab:ParamVal}). It was refined by performing the algorithm on LH time series characterized by a pulsatile pattern with an asymmetric shape of pulses and some regularity in the time evolution of the pulse frequency. LH can be considered as the paragon of any hormone whose secretion pattern is pulsatile, so that the algorithm would also be suited for other hormones (e.g. insulin or growth hormone). As for LH, one has to go through the whole steps of the algorithm, including the removing of 3-point peak (needed parameter: $\lambda_{3p}$) and tunnel-based identification of outliers (needed parameters:  $\alpha$ and $\beta$) for GH and insulin series analysis. In the case of time series for which there is no underlying assumption of asymmetric shape of the peaks or frequency regularity, such as intracellular calcium series, one only needs to go through the first to the fourth steps of the algorithm.

On a more theoretical ground, an interesting question may be addressed in relation to the discretization of a continuous signal. A time series results from a sampling process applied to a continuous signal, which implies that we have chosen (by default) to retrieve the sampling time corresponding to a local maximum to define the time of pulse occurrence. Thus, each IPI is a multiple of the sampling period. As illustrated in the first section, the corresponding theoretical pulse time differs from the pulse occurrence. A deeper analysis of the effect of the sampling on the pulse shape could be undertaken.
This problem is hard to tackle since it mixes non linear dynamics, stochastic process and statistical inference. However, results on this subject would give precious additional knowledge on the location of the theoretical pulse time and could provide more acurate information on the IPI sequence and the frequency encoding.

\section*{Acknowledgments}

The authors wish to thank Dr. Alain Caraty for providing them with additional LH time series. \\
This work is part of the large-scale initiative REGATE (REgulation of the GonAdoTropE axis): \\
\href{http://www.rocq.inria.fr/sisyphe/reglo/regate.html}{http://www.rocq.inria.fr/sisyphe/reglo/regate.html}.

\end{document}